\documentclass[a4paper]{article}
\usepackage{geometry}
\usepackage[colorlinks=true,citecolor=blue]{hyperref}
\usepackage[dvips]{epsfig}
\usepackage{color}
\usepackage{bm}
\usepackage{amssymb}
\usepackage{amsmath}
\usepackage{amsfonts}
\usepackage{subfigure}
\usepackage{xspace}
\usepackage{upgreek}
\usepackage{graphicx}
\usepackage{epstopdf}
\usepackage{float}
\usepackage[table]{xcolor}

\title{Non-polynomial ENO and WENO finite volume methods for hyperbolic conservation laws} 

\author{Jingyang Guo\thanks{Department of Mathematics, University at Buffalo, SUNY, Buffalo, NY 14260-2900, USA. ({\tt jguo4@buffalo.edu}).} \and 
            Jae-Hun Jung\thanks{Corresponding author. Department of Mathematics, University at Buffalo, SUNY, Buffalo, NY 14260-2900, USA. ({\tt jaehun@buffalo.edu}).} }

\begin{document}
\maketitle

\section*{Abstract}
The essentially non-oscillatory (ENO) method is an efficient high order numerical method for solving hyperbolic conservation laws designed to reduce the Gibbs oscillations, if existent, by adaptively choosing the local stencil for the interpolation. The original ENO method is constructed based on the polynomial interpolation and the overall rate of convergence provided by the method is uniquely determined by the total number of interpolation points involved for the approximation. In this paper, we propose simple non-polynomial ENO and weighted ENO (WENO) finite volume methods in order to enhance the local accuracy and convergence. We first adopt the infinitely smooth radial basis functions (RBFs) for a non-polynomial interpolation. Particularly we use the multi-quadric and Gaussian RBFs. The non-polynomial interpolation such as the RBF interpolation offers the flexibility to control the local error by optimizing the free parameter. Then we show that the non-polynomial interpolation can be represented as a perturbation of the polynomial interpolation. That is, it is not necessary to know the exact form of the non-polynomial basis for the interpolation.  In this paper, we formulate the ENO and WENO methods based on the non-polynomial interpolation and derive the optimization condition of the perturbation. To guarantee the essentially non-oscillatory property, we switch the non-polynomial reconstruction to the polynomial reconstruction adaptively near the non-smooth area by using the monotone polynomial interpolation method. The numerical results show that the developed non-polynomial ENO and WENO methods enhance the local accuracy. 

\textbf{keywords}
Essentially non-oscillatory method, Weighted essentially non-oscillatory method, Radial basis function interpolation, Finite volume method, Hyperbolic conservation laws. 


\pagestyle{myheadings}
\thispagestyle{plain}
\markboth{JINGYANG GUO \& JAE-HUN JUNG}{NON-POLYNOMIAL ENO \& WENO METHODS }

\newtheorem{thm}{Theorem}
\newtheorem{example}{Example}
\newtheorem{algorithm}{Algorithm}
\newtheorem{remark}{Remark}

\section{Introduction}
We consider the hyperbolic conservation laws
\begin{eqnarray}
  u_t  + \triangledown \cdot F(u) = 0,  \label{eq1}
\end{eqnarray}
with the state vector u $\equiv$ u(t,x) : $\emph{I} \times \Omega \rightarrow \mathbb{R}^{m}$, for a time interval $\emph{I} :=(0,T]$ for some $T>0$ and an open bounded computational domain $\Omega \subset \mathbb{R}^d$. $F(u) := (f_1(u), \cdots , f_m(u))$ is the flux function. An initial condition $u_0(x) = u(x,0)$ is given along with appropriate boundary conditions. Despite the smoothness of $u_0(x)$, the solution to (\ref{eq1}) may develop a discontinuity within a finite time. High order numerical approximations of the developed discontinuity suffer from the Gibbs phenomenon yielding spurious oscillations near the discontinuity. These oscillations not only degrade the accuracy of the approximation but may also cause instability of the scheme. Since the publications by Harten et al. \cite{Harten} and by Jiang and Shu \cite{WENO}, the essentially non-oscillatory (ENO) and weighted essentially non-oscillatory (WENO) methods have been one of the most powerful numerical methods that can successfully deal with the Gibbs oscillations and widely used in various applications. Numerous modifications of the original ENO and WENO methods have been also developed. These include recent works such as WENO-M \cite{WENOM}, WENO-Z \cite{WENOZ}, power-ENO \cite{Serna}, WENO-P \cite{Ha}, modification of the ENO basis \cite{Berzins} and WENO-$\eta$ methods \cite{Fan}, to name a few. There is no {\it best} ENO/WENO variation because all variations have their own strengths and weaknesses. However, most variations have the common element for their construction: the polynomial reconstruction.  In recent reviews of the WENO methods by Shu \cite{ShuSIAM}, the ENO/WENO reconstruction based on the non-polynomial functions such as the Fourier functions is briefly mentioned such as the one used in \cite{Christofi}. However, most ENO/WENO variations are based on the polynomial reconstruction. In this paper, let us call the ENO/WENO method based on the polynomial reconstruction as the regular ENO/WENO method. 

In this paper, we present a simple new type of the ENO/WENO method based on the non-polynomial interpolation.  As an example of the non-polynomial bases, we first use radial basis functions (RBFs). As we will show later in this paper, the choice of RBFs as a non-polynomial basis function is not necessary but the perturbation form works in more general sense. Furthermore, the presented method is not a hybrid method that combines the ENO/WENO method with the RBF method or other high order methods. The presented method is basically the ENO/WENO method but with the interpolation coefficients slightly modified. In the seminal work of \cite{RBF-ENO}, the ADER method was developed based on the polyharmonic spline, which belongs to the family of piecewise smooth RBFs. The motivation of the method presented in \cite{RBF-ENO} was to adopt the WENO method efficiently for the arbitrary geometry and unstructured mesh by using the meshless feature of RBFs. Our main motivation in this paper, however, is not in using the meshless feature of RBFs but is to enhance the original ENO/WENO accuracy by modifying the interpolation coefficients. For this, we need free parameters and need to optimize them. In \cite{RBF-ENO} the polyharmonic spline was used and there was no undetermined shape parameter --- or the shape parameter is fixed as $\epsilon = 1$ and the order of convergence is overall fixed once the size of each stencil $k$ is fixed. But in our proposed method, the shape parameter is essential. Thus the presented method in this paper is different from the one in  \cite{RBF-ENO}.

\begin{table}[h]
\renewcommand{\arraystretch}{1.5}
\caption{Commonly used radial basis functions $\phi(r), ~ r \geqslant 0$ with $\epsilon$ known as the shape parameter.}
\begin{center} \footnotesize
\begin{tabular}{|c|c|c|l|} \hline  
\multicolumn{2}{|c|}{\textbf{Infinitely smooth RBFs}} & \multicolumn{2}{|c|}{\textbf{Piecewise smooth RBFs}}\\
\hline
Gaussian (GA) & $e^{-(\epsilon r)^2}$ &  Polyharmonic spline & $r^k, ~ k=1,3,5,...$\\ 
Multiquadratic (MQ) & $\sqrt{1+(\epsilon r)^2}$ &   & $r^{k} \ln(r), ~ k=2,4,6,... $ \\  
Inverse quadratic (IQ) & $\frac{1}{1+(\epsilon r)^2}$ &  &  \\
\hline
\end{tabular}
\end{center} 
\label{RBFtabel}
\end{table}

RBFs are divided into two categories depending on whether there are undetermined shape parameters: piecewise smooth RBFs and infinitely smooth RBFs (see Table \ref{RBFtabel}). In this paper, we first use the infinitely smooth RBFs because they are defined with a free parameter $\epsilon$, so-called the shape parameter. Since the parameter is free yet to determined locally it yields the flexibility to improve the original ENO/WENO accuracy. 
%
%
%
We will show the equivalence of the derived interpolations by different RBFs. This means that irrespective of the bases used, we end up with the same kind of reconstruction results. This is also true for the piecewise smooth RBF basis, if we regard them as a special case of the infinitely smooth RBF basis with the shape parameter fixed as $\epsilon = 1$. We then show that the derived RBF interpolation formulas are equivalent to the perturbed polynomial interpolation. Thus one can use other non-polynomial bases rather than RBFs as long as the new basis is defined with one or more free parameters for improving the local accuracy and convergence. For the RBF interpolation, it becomes a polynomial interpolation if we set the shape parameter to vanish. This makes it easy to modify the existing ENO/WENO code to the proposed ENO/WENO method. We restrict our discussion to the one-parameter perturbation although it is possible to utilize multiple free  parameters. The local shape parameter values are determined in such a way that the leading error term in the Taylor series of the reconstruction around each node vanishes. 
%
%
%

Unlike the polynomial interpolation, the perturbed interpolation such as the RBF interpolation is not necessarily  {\it consistent}. Thus it may yield oscillatory interpolation even for constant functions. Such an inconsistency helps the proposed method to enhance local accuracy if the solution is smooth. However, if the solution contains discontinuities the inconsistent reconstruction causes the Gibbs oscillations. To prevent the Gibbs oscillations, we adopt the monotone interpolation method by measuring the local extrema. If the local extrema exists within the interpolating cells, we switch the non-polynomial interpolation into the polynomial interpolation. The switch can be done efficiently by adopting the vanishing shape parameter which reduces the method into a polynomial method \cite{Larson}. Thus, by making the shape parameter vanish we can easily switch the RBF-ENO method to the ENO method and let the regular ENO reconstruction prevent oscillations. 

The paper is composed of the following sections. 
In Section 2, we briefly explain the finite volume ENO and WENO methods.  In Section 3, we introduce the RBF-ENO interpolation based on the multi-quadric (MQ) and Gaussian RBFs and the perturbed interpolation for $k = 2$ and $k = 3$. In this section, we provide the table of the interpolation coefficients for $k = 2$ and $k = 3$. In Section 4, we explain the monotone polynomial interpolation method. The monotone polynomial interpolation is constructed in order to use the vanishing condition of the shape parameter, for which the vanishing condition is derived.
In Section 5, we briefly explain the time-integration and flux schemes that we used for the numerical experiment. Then we present the 1D numerical examples for both linear and nonlinear problems and for scalar and system problems. In Section 6, we explain the 2D ENO/WENO finite volume interpolation method based on the non-polynomial bases. In Section 7 the 2D numerical examples are presented. In Section 8, we provide a brief conclusion and our future research.

\section{Finite volume ENO/WENO method}
Given a grid with $N$ number of points
$$ a = x_{\frac{1}{2}} < x_{\frac{3}{2}} < \cdots < x_{N-\frac{1}{2}} < x_{N+ \frac{1}{2}} = b,$$
for the $i$-th cell $I_i = [x_{i-\frac{1}{2}}, x_{i+\frac{1}{2}}]$, define the cell center $x_i$ and cell size $\Delta x_i$  as
$$x_i = \frac{1}{2} (x_{i-\frac{1}{2}}+ x_{i+\frac{1}{2}}),~  \Delta x_i = x_{i+\frac{1}{2}} - x_{i-\frac{1}{2}},~  i = 1,2,\cdots,N. $$ 
Denote the maximum cell size by
$$ \Delta x =\max_{1 \leq i \leq N} \Delta x_i. $$
At each $i$-th cell, the cell average ${\bar v}_i$ of a function $v(x)$ is given as 
$$ {\bar v}_i = {1\over {\Delta x_i}} \int^{x_{i+{1\over 2}}}_{x_{i - {1\over 2}}} v(\xi) d\xi, \quad i = 1, \cdots, N.$$
For the finite volume ENO method, we seek a function $p_i(x)$ such that we have $k$-th order or higher order accurate approximation to the function $v(x)$ in $I_i$, that is, 
$$
  p_i(x) = v(x) + O(\Delta x^m), \quad {m\ge k}, \quad x \in I_i. 
$$
Note that $m = k$ for the regular ENO method. In this paper we want to find an approximation that yields $m \ge k$ for the proposed method. 
The cell boundary values of $v(x)$ in $I_i$ are then approximated by $p_i(x)$ as 
$$
       v^-_{i+{1\over 2}} = p_i (x_{i+{1\over 2}}) \quad \mbox{ and } \quad  v^+_{i-{1\over 2}} = p_i(x_{i-{1\over 2}}), 
$$
so that they are at least $k$-th order accurate. Here the superscripts $+$ and $-$ denote the right hand side and left hand side limits. 

For the $k$-th order ENO reconstruction, we choose the stencil based on $r$ cells to the left and $s$ cells to the right including $I_i$ such that 
$$
   r + s + 1 = k. 
$$
Define $S_r(i)$ as the stencil composed of those $k$ cells including the cell $I_i$
\begin{equation}
      S_r(i) = \left\{ I_{i-r}, \cdots, I_{i+s} \right\}, ~\quad r=0, \cdots, k-1.
      \label{stencil}
\end{equation}
Define a primitive function $V(x)$ such that 
\begin{equation}
      V(x) = \int^x_{x_{i-r-{1\over 2}}} v(\xi) d\xi, 
      \label{primitive}
\end{equation}
where the lower limit in the integral can be any cell boundary \cite{Shu}. By the definition of $V(x)$ in \eqref{primitive}, $V'(x) = v(x)$ \cite{Shu}. Then for $i - r - 1 \leqslant l \leqslant i + s$, $V(x_{l+{1\over 2}})$ is given by the linear sum of cell averages 
$$
   V(x_{l + {1\over 2}}) = \sum_{j = i-r}^l\int^{x_{j+{1\over 2}}}_{x_{j - {1\over 2}}} v(\xi)d\xi = \sum_{j = i-r}^l \Delta x_j {\bar v}_j  = \sum_{j = i-r}^l \Delta x {\bar v}_j. 
$$
The regular ENO method constructs the polynomial interpolation of $V(x)$ based on $V(x_{l + {1\over 2}}),~ i - r - 1 \leqslant l \leqslant i + s$, while the non-polynomial ENO method constructs the non-polynomial interpolation of $V(x)$ such as the RBF interpolation of $V(x)$.
Suppose $P(x)$ is some interpolation for $V(x)$ such that 
\begin{equation}
      P(x) =  V(x) + O(\Delta x^{m+1}),
      \label{Px}
\end{equation}
then $p(x) \equiv P'(x)$ is the function we seek to approximate $v(x)$ where 
\begin{equation}
      p(x) =  v(x) + O(\Delta x^{m}).
      \label{px}
\end{equation}
We will show that whether we use the polynomial interpolation or the non-polynomial interpolation such as the RBF interpolation, the reconstruction is given in the same form as below 
$$v^{(r)-}_{i+ {1\over 2}} \equiv p(x_{i+{1\over 2}}) = \sum^{k-1}_{j=0} c_{rj} {\bar{v}}_{i-r+j}$$ 
with slightly different reconstruction coefficients $c_{rj}$.

The WENO reconstruction is then formulated as a convex combination of all possible ENO reconstruction in $S_r(i)$. In $S_r(i)$ \eqref{stencil} there are $k$ different ENO reconstructions of $v^{(r)-}_{i+{1\over 2}}$ and the WENO reconstruction would take the convex combination of all those reconstructions of $v^{(r)-}_{i+ {1\over 2}}$:
\begin{eqnarray}
   v^-_{i+{1\over 2}} = \sum^{k-1}_{r=0} w_r v^{(r)-}_{i+ {1\over 2}},
\label{weno_re}
\end{eqnarray}
where 
$$
   w_r = \frac{\alpha_r}{\sum^{k-1}_{s=0} \alpha_s}, r=0,\cdots ,k-1
$$
with 
$$
   \alpha_r = \frac{d_r}{(\epsilon+\beta_r)^2}.
$$
Here $d_r$ are the polynomial interpolation coefficients and $\epsilon>0$ is introduced to avoid the case that the denominator becomes $0$ which is usually taken as $10^{-6}$. $\beta_r$ are the ``smooth indicators" of the stencil $S_r(i)$. For the details of $d_r$ and $\beta_r$, we refer readers to \cite{Shu}.

\section{1D interpolation} 
We first consider the case with $k=2$, that is, two cells are used for the ENO reconstruction. For this case, three cell averages ${\bar v}_{i-1}$, ${\bar v}_{i}$ and ${\bar v}_{i+1}$ are available. To reconstruct the boundary values of $v^-_{i+ \frac{1}{2}}$ and $v^+_{i- \frac{1}{2}}$, we use either $\{ {\bar v}_{i-1}, {\bar v}_{i} \}$ or $\{ {\bar v}_{i}, {\bar v}_{i+1}\}$. Which cell averages should be used is decided by the Newton's divided difference method \cite{Shu}. Suppose that for the given cell $I_i$ the Newton's divided difference method suggests that we use $\{ {\bar v}_{i}, {\bar v}_{i+1}\}$ for the local reconstruction from the available stencil $\{ {\bar v}_{i-1}, {\bar v}_{i}, {\bar v}_{i+1} \}$. For simplicity, we use the uniform grid, i.e. $ \Delta x_i = \Delta x, \forall i$. We only show the reconstruction at the cell boundary $x = x_{i+ {1\over 2}}$. The reconstruction for the other cell boundary $x = x_{i- {1\over 2}}$ can be obtained in the same manner.

Define $V(x) = \int^x_{x_{i+ {1\over 2}}} v(\xi) d\xi $. Then the primitive function at the cell boundaries are given by 
\begin{eqnarray}
V(x_{i-{1\over 2}}) &=& - \Delta x \cdot {\bar v}_{i} \nonumber\\
V(x_{i+{1\over 2}}) &=& 0, \nonumber \\
V(x_{i+{3\over 2}}) &=& \Delta x \cdot {\bar v}_{i+1}. \nonumber 
\end{eqnarray}
\subsection{Polynomial reconstruction}
The polynomial interpolation $P(x)$ of $V(x)$  is given as 
$$ P(x) = \lambda_1 + \lambda_2 x + \lambda_3 x^2. $$
Let $ {\vec{V}} =  [V_{i- {1\over 2}}, V_{i+ {1\over 2}}, V_{i+ {3\over 2}} ]^T$, $ {\vec{\lambda}} = [\lambda_1, \lambda_2,  \lambda_3]^T$ and the interpolation matrix $A$ be
$$
A = \begin{bmatrix} 
1 & x_{i-{1\over 2}} & x_{i-{1\over 2}}^2
\\ 
1 & x_{i+{1\over 2}} & x_{i+{1\over 2}}^2
\\ 
1 & x_{i+{3\over 2}} & x_{i+{3\over 2}}^2
\end{bmatrix}. 
$$
Then the expansion coefficients $\lambda_i$ are determined by solving the linear system $ \vec{V} = A \cdot \vec{\lambda} $.
After taking the first derivative of $P(x)$ and plugging $x = x_{i+ {1\over 2}}$, we obtain
\begin{equation}
v^-_{i+ {1\over 2}} = p(x_{i+ {1\over 2}}) = {1\over 2} \cdot \bar{v}_i+{1\over 2} \cdot \bar{v}_{i+1}. 
\label{poly_recon}
\end{equation}
The Taylor series expansion of $v^-_{i+ {1\over 2}}$ around $x=x_{i+{1\over 2}}$ yields
\begin{equation}
v^-_{i+ {1\over 2}} = p(x_{i+{1\over 2}}) = v_{i+{1\over 2}} + {1\over 6}  v''_{i+{1\over 2}} \Delta x^2 + {1\over 120} v^{(4)}_{i+{1\over 2}} {\Delta x^4} + O(\Delta x^6). 
\label{poly_recon_error}
\end{equation}
The first term in the right hand side in \eqref{poly_recon_error} is the exact value of $v(x)$ at $x = x_{i+ {1\over 2}}$. So we confirm that \eqref{poly_recon} is a 2nd order reconstruction.

\subsection{Multiquadratic RBF reconstruction}

For a non-polynomial reconstruction, we consider the MQ-RBF interpolation. The MQ-RBF interpolation $P(x)$ for $V(x)$ is given by
$$ P(x) = \lambda_1 \sqrt{1+\epsilon^2(x-x_{i- {1\over 2}})^2  } + \lambda_2 \sqrt{1+\epsilon^2(x-x_{i+ {1\over 2}})^2  } + \lambda_3 \sqrt{1+\epsilon^2(x-x_{i+ {3\over 2}})^2  }.  $$
Then the interpolation matrix $A$ is given by 
$$
A = \begin{bmatrix} 
1 & \sqrt{\Delta x^2 \epsilon_2^2+1} & \sqrt{4 \Delta x^2 \epsilon_3^2+1}
\\ 
\sqrt{\Delta x^2 \epsilon_1^2+1} & 1 & \sqrt{\Delta x^2 \epsilon_3^2+1}
\\ 
\sqrt{4 \Delta x^2 \epsilon_1^2+1} & \sqrt{\Delta  x^2 \epsilon_2^2+1} & 1
\end{bmatrix}. 
$$
Again, we take the first derivative of $P(x)$ and plug $x=x_{i+{1\over 2}}$ to obtain
\begin{equation}
v^-_{i+ {1\over 2}} = p(x_{i+ {1\over 2}}) = \frac{\sqrt{4 \epsilon^2 \Delta x^2  +1}+1}{4 \sqrt{\epsilon^2 \Delta x^2 +1} } \cdot \bar{v}_i+\frac{\sqrt{4 \epsilon^2 \Delta x^2 +1}+1}{4 \sqrt{\epsilon^2 \Delta x^2 +1} } \cdot \bar{v}_{i+1} , 
\label{i+1/2_exact}
\end{equation}
Expanding $v^-_{i+ {1\over 2}}$ around $x=x_{i+{1\over 2}}$ in the Taylor series yields
\begin{eqnarray}
v^-_{i+ {1\over 2}} = p(x_{i+{1\over 2}}) &=& v_{i+{1\over 2}} + \left({1\over 2}\epsilon^2 v_{i+{1\over 2}}+{1\over 6}  v''_{i+{1\over 2}} \right) \Delta x^2 
\nonumber \\
& & + \left( -{9\over 8} \epsilon^4  v_{i+{1\over 2}}  + {1\over 12} \epsilon^2  v''_{i+{1\over 2}} + {1\over 120} v^{(4)}_{i+{1\over 2}}  \right) {\Delta x^4} + O(\Delta x^6). 
\label{p(i+1/2)for_error}
\end{eqnarray}
Thus we see that \eqref{i+1/2_exact} is at least $2$nd order accurate to $v_{i+{1\over 2}}$. If we take the value of $\epsilon$ as below
\begin{equation}
\epsilon^2 = -{1\over 3}{{v''_{i+{1\over 2}}}\over{v_{i+{1\over 2}}}},
\label{i+1/2_eps}
\end{equation}
then we obtain a $4$th order accurate approximation, i.e. $p(x_{i+{1\over 2}}) = v_{i + {1\over 2}} + O(\Delta x^4)$. 

We notice that the coefficients of ${\bar v}_i$ and ${\bar v}_{i+1}$ in  \eqref{i+1/2_exact} are in a complicated form involving a calculation of square roots. 
We can alternatively drive more efficient and simpler forms that can still yield the same desired order. To do this, we use the Taylor series again to expand the right hand sides of \eqref{i+1/2_exact} as below
\begin{eqnarray}
v^-_{i+ {1\over 2}} = p(x_{i+ {1\over 2}}) &=& \left(\frac{1}{2}+\frac{1}{4} \epsilon^2 \Delta x^2-\frac{9}{16} \epsilon^4 \Delta x^4\right) \cdot \bar{v}_i + \left(\frac{1}{2}+\frac{1}{4} \epsilon^2 \Delta x^2-\frac{9}{16} \epsilon^4 \Delta x^4\right) \cdot \bar{v}_{i+1} \nonumber \\
  &&+ O(\Delta x^6), 
\label{i+1/2_around_dx}
\end{eqnarray}
Since we want the RBF-ENO method with $k=2$ to give $3$rd order accuracy at least, we can ignore all the high order terms in \eqref{i+1/2_around_dx}, which yields 
\begin{equation}
v^-_{i+ {1\over 2}} = p(x_{i+ {1\over 2}}) = \left(\frac{1}{2}+\frac{1}{4} \epsilon^2 \Delta x^2\right) \cdot \bar{v}_i + \left(\frac{1}{2}+\frac{1}{4} \epsilon^2 \Delta x^2\right) \cdot \bar{v}_{i+1}. 
\label{i+1/2_approx}
\end{equation}
Carry out the error analysis on our new reconstructions \eqref{i+1/2_approx} and we get
\begin{eqnarray}
v^-_{i+ {1\over 2}} = p(x_{i+{1\over 2}}) &=& v_{i+{1\over 2}} + \left({1\over 2}\epsilon^2  v_{i+{1\over 2}}  + {1\over 6}  v''_{i+{1\over 2}} \right) \Delta x^2  
\nonumber \\
&& + \left( {1\over 12} \epsilon^2  v''_{i+{1\over 2}} + {1\over 120} v^{(4)}_{i+{1\over 2}} \right) {\Delta x^4} + O(\Delta x^6). 
\label{p(i+1/2)for_error_approx}
\end{eqnarray}
We realize that getting rid of all the high order terms $O(\Delta x^4)$ in \eqref{i+1/2_around_dx} would only exert an influence on $O(\Delta x^4)$ terms in \eqref{p(i+1/2)for_error}, while all the $O(\Delta x^2)$ terms stay the same. In another word, we can still use \eqref{i+1/2_eps} to remove the $2$nd order terms and achieve the $3$rd order accuracy for the new reconstruction. 

Now the problem is how to calculate such values of $\epsilon$ in \eqref{i+1/2_eps} to achieve the higher order accuracy than the $2$nd order. To explain this, we use the point $x = x_{i+{1\over 2}}$ as an example. The case of $x=x_{i-{1\over 2}}$ can be done in the same way. Our idea is as following. Although $v^-_{i+{1\over 2}}$ is computed based on ${\bar v}_{i}$ and ${\bar v}_{i+1}$, we notice that the cell average to the left of $I_i$, i.e. ${\bar v}_{i-1}$ is also known -- notice that it was already used when the Newton's divided difference method was applied to determine the cell averages to be used for the reconstruction for the cell $I_i$. Thus, we can approximate $v''_{i+{1\over 2}}$ and $v_{i+{1\over 2}}$ if all these  given cell average information is used. Construct the Lagrange interpolation of $V(x)$ based on $V_{i-{3\over 2}}, V_{i-{1\over 2}}, V_{i+{1\over 2}},V_{i+{3\over 2}}$ and take the first derivative of $V$ at $x = x_{i+{1\over 2}}$. Then we have 
$$
    v_{i+{1\over 2}} = -{1\over 6} {\bar v}_{i-1} + {5\over 6} {\bar v}_{i} + {1\over 3} {\bar v}_{i+1} + O(\Delta x^3). 
$$
The second derivative of $ v_{i+{1\over 2}} $ is approximated by the third derivative of $V(x)$, which is given by
$$
  v''_{i+{1\over 2}} = {{{\bar v}_{i-1} - 2 {\bar v}_{i} + {\bar v}_{i+1}}\over{\Delta x^2}}+ O(\Delta x). 
$$
Then by plugging the above approximations of $v_{i+{1\over 2}}$ and $v''_{i+{1\over 2}}$ into \eqref{i+1/2_eps}, we determine the value of $\epsilon^2$ as below
\begin{equation}
\epsilon^2 \approx  \frac{2}{\Delta x^2}\cdot \frac{-{\bar v}_{i-1}+2{\bar v}_{i}-{\bar v}_{i+1}}{-{\bar v}_{i-1}+5{\bar v}_{i}+2 {\bar v}_{i+1}+ \epsilon_M}.
\label{i+1/2_eps_approx}
\end{equation}
$\epsilon_M$ is a positive small number to avoid the denominator being zero. Here note that $\epsilon$ can be a complex number because $\epsilon^2$ can be a negative number but it does not harm the RBF interpolation because it is $\epsilon^2$ not $\epsilon$ that is used in the RBF-ENO reconstruction and all the operations are done on real numbers. If we replace $\epsilon^2$ in \eqref{p(i+1/2)for_error} with that in \eqref{i+1/2_eps_approx}, we get the following 
\begin{equation}
 v^-_{i+ {1\over 2}} = p(x_{i+{1\over 2}}) = v_{i+{1\over 2}} + {1\over 12} v^{(3)}_{i+{1\over 2}} \Delta x^3 + O(\Delta x^4). 
\end{equation}
Thus we confirm that the above approximate value of $\epsilon$ in \eqref{i+1/2_eps_approx} actually yields the $3$rd order accuracy.

\subsection{Gaussian RBF reconstruction}
Now we take a look at another infinitely smooth RBF basis, i.e. the Gaussian RBF. 
We want to show that all the RBF reconstructions are equivalent regardless of the basis used.
The Gaussian RBF interpolation $P(x)$ of $V(x)$ is given as
$$ P(x) = \lambda_1 e^{-\epsilon^2(x-x_{i- {1\over 2}})^2  } + \lambda_2 e^{-\epsilon^2(x-x_{i+ {1\over 2}})^2  } + \lambda_3 e^{-\epsilon^2(x-x_{i+ {3\over 2}})^2  }.  $$
Follow the same procedure as before. The exact RBF reconstruction at $x=x_{i+{1\over 2}}$ is
\begin{equation}
v^-_{i+ {1\over 2}} = p(x_{i+ {1\over 2}}) = {\frac{2 \epsilon^2 \Delta x^2 e^{-\epsilon^2 \Delta x^2}}{1-e^{4 \epsilon^2 \Delta x^2}}} \cdot \bar{v}_i+{\frac{2 \epsilon^2 \Delta x^2 e^{-\epsilon^2 \Delta x^2}}{1-e^{4 \epsilon^2 \Delta x^2}} } \cdot \bar{v}_{i+1} . 
\label{i+1/2_exact1}
\end{equation}
Expand the reconstruction in the Taylor series and ignore all the high order terms to get 
\begin{equation}
v^-_{i+ {1\over 2}} = p(x_{i+ {1\over 2}}) = \left(\frac{1}{2}+\frac{1}{2} \epsilon^2 \Delta x^2\right) \cdot \bar{v}_i + \left(\frac{1}{2}+\frac{1}{2} \epsilon^2 \Delta x^2\right) \cdot \bar{v}_{i+1}.
\label{i+1/2_approx1}
\end{equation}
Expanding ${\bar v}_i$ and ${\bar v}_{i+1}$ around the boundary value yields 
\begin{eqnarray}
v^-_{i+ {1\over 2}} = p(x_{i+{1\over 2}}) &=& v_{i+{1\over 2}} + \left(\epsilon^2  v_{i+{1\over 2}}  + {1\over 6}  v''_{i+{1\over 2}} \right) \Delta x^2  
\nonumber \\
&& + \left( {1\over 12} \epsilon^2  v''_{i+{1\over 2}} + {1\over 120} v^{(4)}_{i+{1\over 2}} \right) {\Delta x^4} + O(\Delta x^6). 
\label{p(i+1/2)for_error_approx1}
\end{eqnarray}
Again, if we have
\begin{equation}
\epsilon^2 =  \frac{1}{\Delta x^2}\cdot \frac{-{\bar v}_{i-1}+2{\bar v}_{i}-{\bar v}_{i+1}}{-{\bar v}_{i-1}+5{\bar v}_{i}+2 {\bar v}_{i+1}},
\label{i+1/2_eps_approx1}
\end{equation}
then \eqref{p(i+1/2)for_error_approx1} becomes
\begin{equation}
 v^-_{i+ {1\over 2}} = p(x_{i+{1\over 2}}) = v_{i+{1\over 2}} + {1\over 12} v^{(3)}_{i+{1\over 2}} \Delta x^3 + O(\Delta x^4).
\end{equation}
Thus we see that no matter what RBF basis we use the Taylor series of each reconstruction results in the same expansion after the proper appropriate value of $\epsilon$ is plugged in.

\subsection{Perturbed polynomial reconstruction}
\label{perturbation}
Now let us consider the generalized case by perturbing the polynomial interpolation instead of using any specific RBF basis. By knowing the RBF interpolations above, assume that the unknown non-polynomial interpolation of the $2$nd order or higher is given as the perturbed form of the polynomial interpolation as below
\begin{equation}
v^-_{i+ {1\over 2}} = p(x_{i+ {1\over 2}}) = \left(\frac{1}{2}+a \epsilon^2 \Delta x^2\right) \cdot \bar{v}_i + \left(\frac{1}{2}+b \epsilon^2 \Delta x^2\right) \cdot \bar{v}_{i+1}, 
\label{i+1/2_approx2}
\end{equation}
where $a$ and $b$ are real constants and $\epsilon^2$ is to be determined. Expand the given cell averages around $x_{i+{1\over 2}}$ and we get
\begin{eqnarray}
v^-_{i+ {1\over 2}} = p(x_{i+{1\over 2}}) &=& v_{i+{1\over 2}} + \left[ (a+b) \epsilon^2  v_{i+{1\over 2}}  + {1\over 6}  v''_{i+{1\over 2}} \right] \Delta x^2  
\nonumber \\
&& - {1\over 2} (a-b) \epsilon^2 v'_{i+{1\over 2}} \Delta x^3  + O(\Delta x^4). 
\label{p(i+1/2)for_error_approx2}
\end{eqnarray}
To make the $2$nd order error term vanish, we take the unknown parameter $\epsilon$ as 
\begin{equation}
\epsilon^2 =  \frac{1}{\Delta x^2 (a+b)}\cdot \frac{-{\bar v}_{i-1}+2{\bar v}_{i}-{\bar v}_{i+1}}{-{\bar v}_{i-1}+5{\bar v}_{i}+2 {\bar v}_{i+1}}. 
\label{i+1/2_eps_approx2}
\end{equation}
Then \eqref{p(i+1/2)for_error_approx2} becomes
\begin{equation}
v^-_{i+ {1\over 2}} = p(x_{i+{1\over 2}}) = v_{i+{1\over 2}} + \left( {1\over 12} v^{(3)}_{i+{1\over 2}} + \frac{a-b}{a+b} \frac{v^{'}_{i+{1\over 2}} v^{''}_{i+{1\over 2}}}{v_{i+{1\over 2}}} \right) \Delta x^3 + O(\Delta x^4).
\label{error_result2}
\end{equation}
There are multiple ways to determine the values of the perturbed coefficients $a$ and $b$. For example, 
\begin{eqnarray}
 \label{(1)} 
a+b\neq 0, \\
\label{(2)}
a-b = 0.   
\end{eqnarray}
Condition \eqref{(1)} implies that the non-polynomial reconstruction is inconsistent. And we see that both \eqref{i+1/2_approx} and \eqref{i+1/2_approx1} satisfy \eqref{(1)} and \eqref{(2)}. In this way, we see that the RBF reconstruction irrespective of the basis used is equivalent to the perturbed polynomial reconstruction. That is, it is not necessary to know the exact form of the non-polynomial basis for the interpolation.  

\subsection{Reconstruction coefficients for $k=2$ and $k=3$}
The left table in Table \ref{table_k2} shows the coefficients of the polynomial reconstruction that yield $2$nd order convergence. The right table provides the coefficients with the RBF reconstruction corresponding to the MQ-RBF interpolation \eqref{i+1/2_approx} which yields $3$rd order convergence. The left table is obtained by the limit of the RBF reconstruction as $\epsilon \rightarrow 0$, the polynomial limit. Notice that the RBF reconstruction is the same as the polynomial reconstruction with the  perturbation term added. Table \ref{table_k3} shows the coefficients for $k = 3$. 

\begin{table}[h]
\renewcommand{\arraystretch}{2.5}
\caption{Left: The polynomial reconstruction coefficients, also for the RBF interpolation with $\epsilon \rightarrow 0$. Right: The MQ-RBF reconstruction coefficients  $c_{rj}$ corresponding to \eqref{i+1/2_approx}. $\eta =  \epsilon^2 \Delta x^2$.}
\begin{center} \footnotesize
\begin{tabular}{|c|c|c|c|} \hline  
k & r & j=0 & j=1\\ \hline 
  & -1 &  $\frac{3}{2}$ &  $-\frac{1}{2}$\\  \cline{2-4}
2 & 0 &  $\frac{1}{2} $ & $\frac{1}{2} $ \\  \cline{2-4} 
  & 1 & $-\frac{1}{2}$ &  $\frac{3}{2}$ \\   \cline{1-4} 
\multicolumn{4}{|c|}{$\eta = 0$} \\
\hline
\end{tabular}
\hskip .2in
\begin{tabular}{|c|c|c|c|} \hline  
k & r & j=0 & j=1\\ \hline 
  & -1 & ${3\over 2} - {3\over 2} \eta$ & $- {1\over 2} + {1\over 2} \eta$\\  \cline{2-4}
2 & 0 &  ${1\over 2} + {1\over 4} \eta$ & $ {1\over 2} + {1\over 4} \eta$ \\  \cline{2-4} 
  & 1 & $- {1\over 2} + {1\over 2} \eta$ &  ${3\over 2} - {3\over 2} \eta$ \\
\cline{1-4} 
\multicolumn{4}{|c|}{$\eta = \epsilon^2 \Delta x^2  = \frac{ 2(-{\bar v}_{i-1}+2{\bar v}_{i}-{\bar v}_{i+1})}{-{\bar v}_{i-1}+5{\bar v}_{i}+2 {\bar v}_{i+1}+ \epsilon_M}$} \\
\hline
\end{tabular}
\end{center} 
\label{table_k2}
\end{table}


\begin{table}[h]
\renewcommand{\arraystretch}{2.5}
\caption{Left: The polynomial reconstruction coefficients, also for the RBF interpolation with $\epsilon \rightarrow 0$. Right: The MQ-RBF reconstruction coefficients  $c_{rj}$, $\eta =  \epsilon^2 \Delta x^2$.}
\begin{center} \footnotesize
\begin{tabular}{|c|c|c|c|c|} \hline  
k & r & j=0 & j=1 & j=2\\ \hline 
  & -1 & $\frac{6}{11}$ & $-\frac{7}{6}$ & $\frac{1}{3}$ \\  \cline{2-5}
3 & 0 &  $\frac{1}{3}$ & $\frac{5}{6}$ & $-\frac{1}{6}$ \\  \cline{2-5} 
  & 1 & $-\frac{1}{6}$ & $\frac{5}{6}$ & $\frac{1}{3}$\\  \cline{2-5}
  & 2 & $\frac{1}{3}$ & $-\frac{7}{6}$ & $\frac{11}{6}$ \\  \cline{1-5} 
\multicolumn{5}{|c|}{$\eta = 0$} \\
\hline
\end{tabular} 
\hskip .2in
\begin{tabular}{|c|c|c|c|c|} \hline  
k & r & j=0 & j=1 & j=2\\ \hline
  & -1 & $\frac{6}{11} - \frac{9}{2} \eta$ & $-\frac{7}{6} + 6 \eta$ & $\frac{1}{3} - \frac{3}{2} \eta$ \\  \cline{2-5}
3 & 0 &  $\frac{1}{3} + \frac{5}{6} \eta$ & $\frac{5}{6} - \frac{2}{3} \eta$ & $-\frac{1}{6} - \frac{1}{6} \eta$ \\  \cline{2-5} 
  & 1 & $-\frac{1}{6} - \frac{1}{6} \eta$ & $\frac{5}{6} - \frac{2}{3} \eta$ & $\frac{1}{3} + \frac{5}{6} \eta$\\  \cline{2-5}
  & 2 & $\frac{1}{3} - \frac{3}{2} \eta$ & $-\frac{7}{6} + 6 \eta$ & $\frac{11}{6} - \frac{9}{2} \eta$ \\  \cline{1-5} 
\multicolumn{5}{|c|}{$\eta =\epsilon^2 \Delta x^2=  {{{\bar{v}}_{i-1} - 3{\bar{v}}_i + 3 {\bar{v}}_{i+1} - {\bar{v}}_{i+2}}\over{{\bar{v}}_{i-1} - 15 {\bar{v}}_i + 15 {\bar{v}}_{i+1} - {\bar{v}}_{i+2} + \epsilon_M} } $} \\
\hline
\end{tabular} 

\end{center} 
\label{table_k3}
\end{table}

\section{Switcing condition of $\epsilon$}
The regular ENO method helps to avoid the interpolation using the cell averages across the discontinuity by adaptively choosing the proper cells using the Newton's divided difference method. This essentially yields the non-oscillatory reconstruction near the discontinuity as $N \rightarrow \infty$. Although the non-polynomial ENO reconstruction also utilizes the Newton's divided difference method to determine the proper cells, the value of $\epsilon^2$ used in the interpolation coefficients is determined by the cell averages across the discontinuity as the WENO reconstruction. This is why the non-polynomial ENO reconstruction achieves higher order accuracy in the smooth area. However, it may not necessarily be non-oscillatory near the discontinuity for the same reason. 

One way to deal with this problem is to make the RBF-ENO reconstruction become the regular ENO reconstruction in the non-smooth area. This can be done by adopting the vanishing shape parameter  $\epsilon^2$ in the non-smooth area. That is, in the non-smooth area, we want to recover at least the ENO reconstruction, i.e. $\epsilon \rightarrow 0$. To achieve this reduction, we first need to identify which cell contains the discontinuity, which requires extra work such as the edge detection algorithm. The edge detection algorithm has been utilized for the construction of hybrid methods such as the spectral-WENO method \cite{BrunoDon,CDGS2004} and the Fourier continuation WENO method \cite{Hesthaven_FC}. For those hybrid methods, the edge detection algorithm identifies the edge location and the neighboring stencils containing the edge cells are treated by the WENO method and the other cells are treated by the spectral or Fourier methods. For example, the edge detection algorithm used in the spectral-WENO hybrid method and the Fourier continuation method is the Harten's multi-resolution analysis \cite{HartenMR}. We could adopt such a method but in this paper we use a simpler monotone polynomial method described below which yet yields our desired results. 

To illustrate the monotone polynomial method, consider a stencil of $3$ cells for $k = 2$. For the RBF-ENO reconstruction, the shape parameter involves the $2$nd derivative of $v$, i.e. $v''$ as obtained in \eqref{i+1/2_eps}
$$
  \epsilon^2 = -{1\over 3}{{  v''_{i+{1\over 2}}}\over{v_{i+{1\over 2}}}}.
$$
If the function $v(x)$ is discontinuous, the value of $\epsilon^2$ is determined by the cell average values across the discontinuity and the second derivative $v''$ of the reconstruction becomes large resulting in an oscillatory reconstruction. Consider the whole stencil $S = \left\{ I_{i-1}, I_i, I_{i+1}\right\}$. The whole interval of $x$ in $S$ is 
$$
 x_{i-{3\over 2}} \le x \le x_{i+{3\over 2}}.
$$
In this interval,  $p(x)$ by the polynomial reconstruction is given by the first derivative of $P(x)$ as the second order polynomial for $k = 2$,
$$
p(x) \in \mathcal{P}_2 = span\{1, x, x^2 \}, 
$$
where $\mathcal{P}_2$ is a set composed of all polynomials of degree at most $2$. 
Without loss of generality, let the grid points be 
\begin{eqnarray}
  x_{i-{3\over 2}} = 0, \quad 
   x_{i+{3\over 2}} = 3 \Delta x. \nonumber 
\end{eqnarray}
Then, the local maximum or minimum of $p(x)$ exists at $x = x_p$
\begin{eqnarray}
x_p = {{-2{\bar v}_{i-1} +3{\bar v}_{i}-{\bar v}_{i+1} }\over{-{\bar v}_{i-1} +2{\bar v}_{i}-{\bar v}_{i+1}}} \Delta x.
\label{root_value}
\end{eqnarray}
If $x_p$ exists inside the given stencil $S$, i.e. 
$$
     0= x_{i-{3\over 2}} < x_p < x_{i+{3\over 2}} = 3 \Delta x, 
$$
or
$$
 {0} < {{-2{\bar v}_{i-1} +3{\bar v}_{i}-{\bar v}_{i+1} }\over{-{\bar v}_{i-1} +2{\bar v}_{i}-{\bar v}_{i+1}}}  < {3},
$$
then $p(x)$ is not a monotone function in $S$. If $x_p$ exists outside the interval, then $p(x)$ is either monotone increasing or monotone decreasing. We make the non-polynomial ENO method reduced to the regular ENO method in the region where $p(x)$ is not monotone. Using this property we apply the following condition for the non-smoothness area
\begin{eqnarray}
  \epsilon^2 = \left\{ \begin{array}{ll} 0  & \mbox{ if }  {0} < {{-2{\bar v}_{i-1} +3{\bar v}_{i}-{\bar v}_{i+1} }\over{-{\bar v}_{i-1} +2{\bar v}_{i}-{\bar v}_{i+1}}}  < {3} \\
   -{1\over 3}{{  v''_{i+{1\over 2}}} / {v_{i+{1\over 2}}}} & \mbox{ otherwise } \end{array} \right. 
\label{adaptation}
\end{eqnarray}
Note that this algorithm not only makes the non-polynomial ENO method reduced to the regular ENO method in the non-smooth area, but also near the critical point. It is possible that we do not achieve the desired higher order accuracy around the critical point with the monotone polynomial method because with this method the ENO reconstruction is recovered around the critical point as well. In \cite{GuoJungJCP2015}, we show that the monotone polynomial method can be modified so that the reduction near the critical point can be avoided. We do not include the results here because it is beyond the scope of this paper.

\section{1D Numerical experiments}
\subsection{Time integration and flux scheme}
Now we apply the developed method to hyperbolic conservation laws \eqref{eq1}. For the cell averaged value on the given cell $I_i$, ${\bar u}_i$, we seek the solution by solving the following numerically
\begin{eqnarray}
  { {d{\bar u}_i} \over {dt} } = -{{f(u(x_{i+{1\over 2}})) - f(u(x_{i-{1\over 2}}))}\over{\Delta x}}.
\label{PDE1}
\end{eqnarray}
For the numerical solution, the exact flux $f$ in \eqref{PDE1} is replaced by the numerical flux function $h$ and the boundary values are computed by the non-polynomial ENO reconstruction or non-polynomial WENO reconstruction. For the flux, we use the Lax-Friedrichs flux scheme
$$
  h(a,b) = {1\over 2} [f(a) + f(b) - \alpha (b - a)], 
$$
where $\alpha = \max_{u} |{{\partial f}\over{\partial u}} | $ is a constant. For the time integration, we use the TVD RK-3 method \cite{RK3}. 

\subsection{Numerical results}
\begin{figure}[h]
\begin{center}
\includegraphics[width=0.47\textwidth]{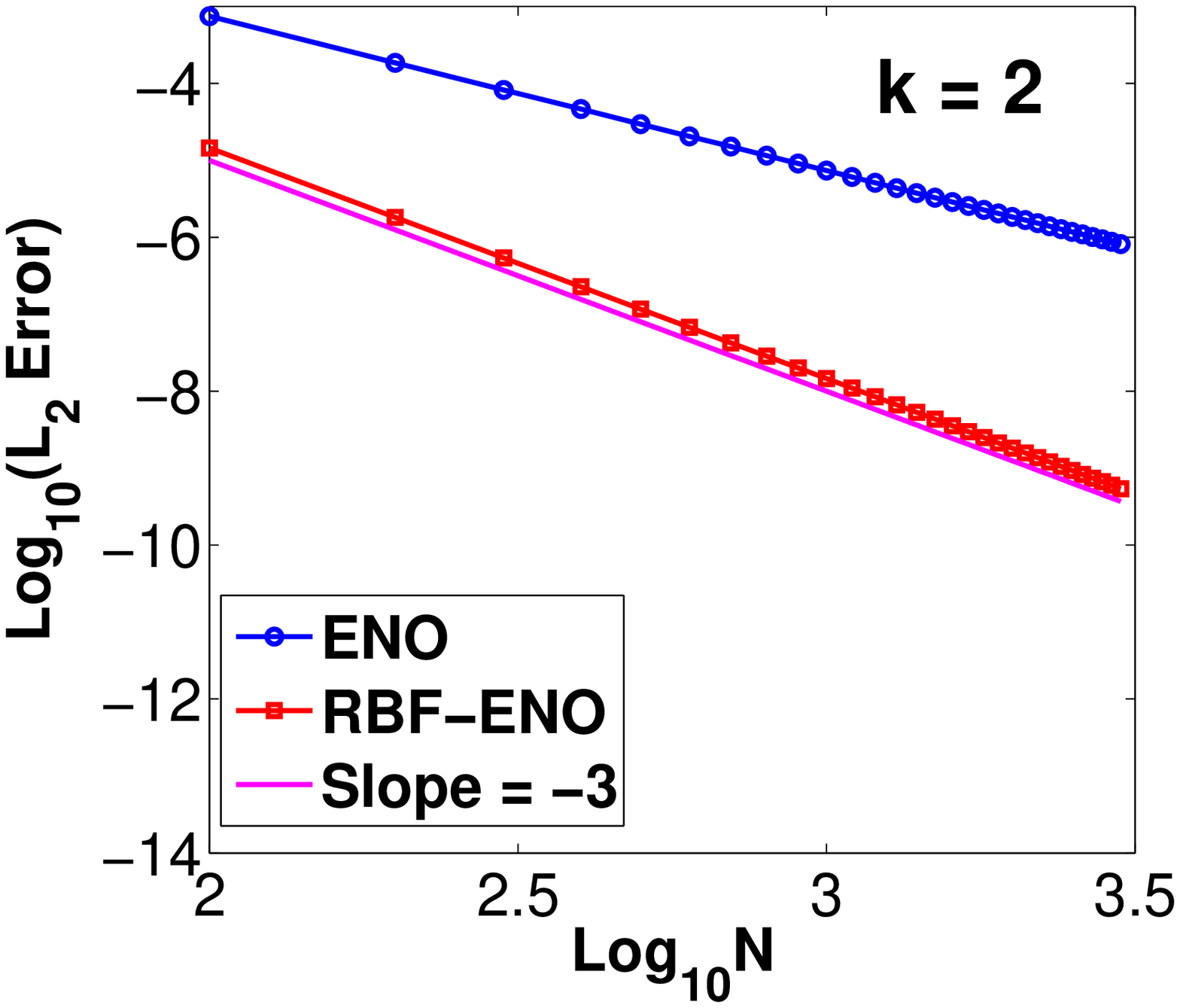}
\includegraphics[width=0.47\textwidth]{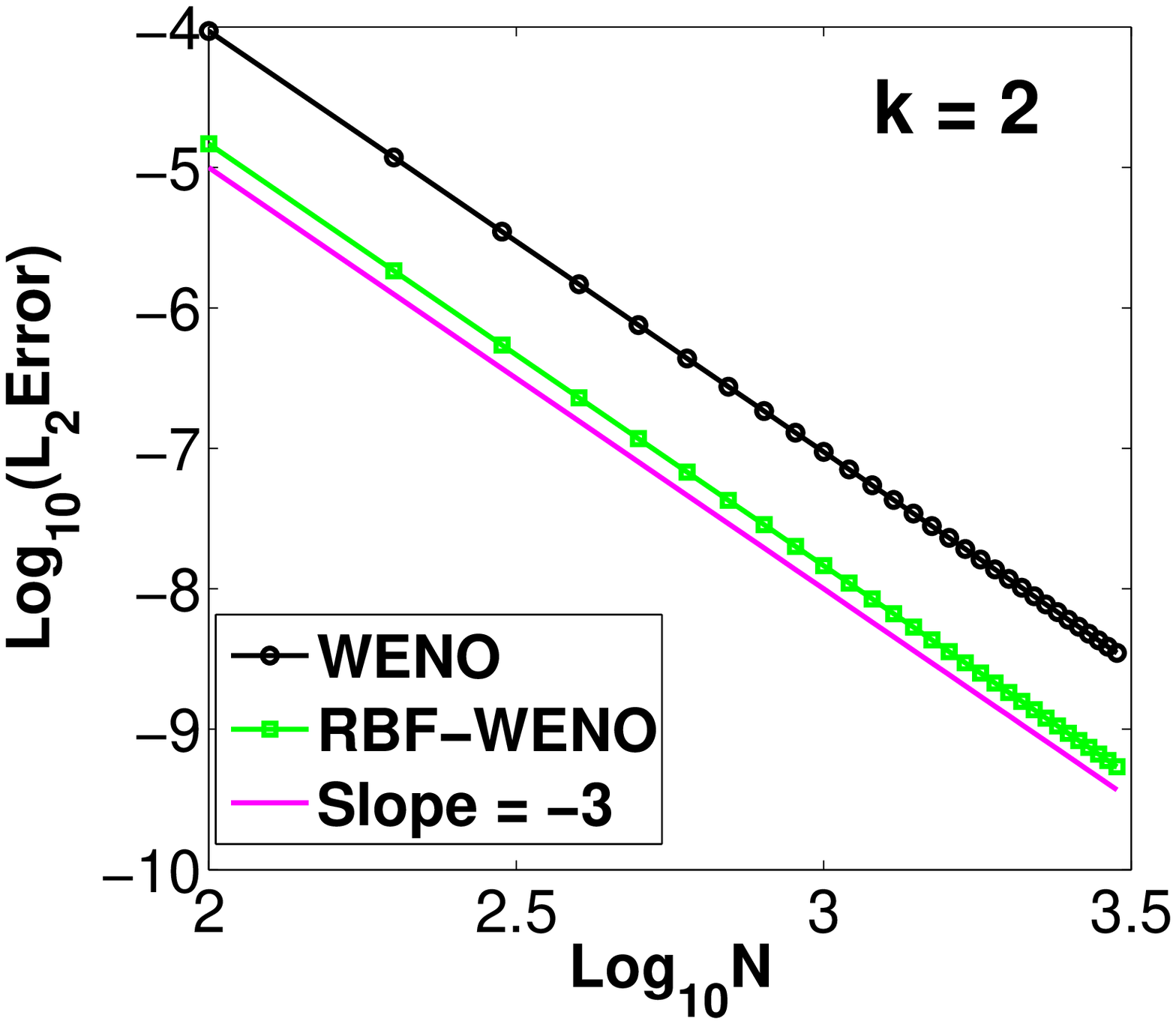}
\end{center}
\caption{(Color online).  $L_2$ errors versus $N$ on logarithmic scale with $k = 2$. Left: ENO method (blue) and RBF-ENO method (red). Right: WENO method (black) and RBF-WENO method (green).}
\label{test1}
\end{figure}

\begin{figure}[h]
\begin{center}
\includegraphics[width=0.47\textwidth]{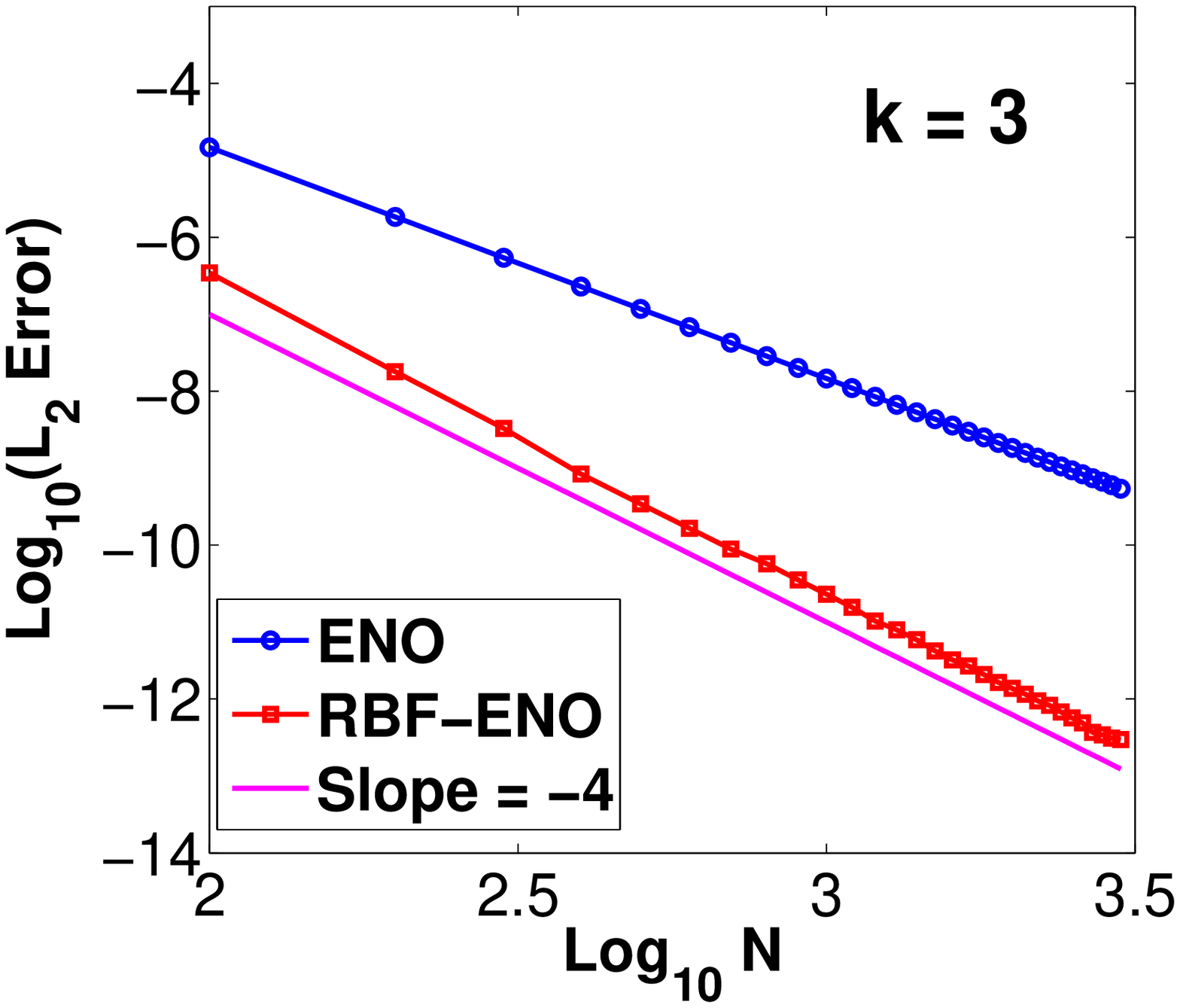}
\includegraphics[width=0.47\textwidth]{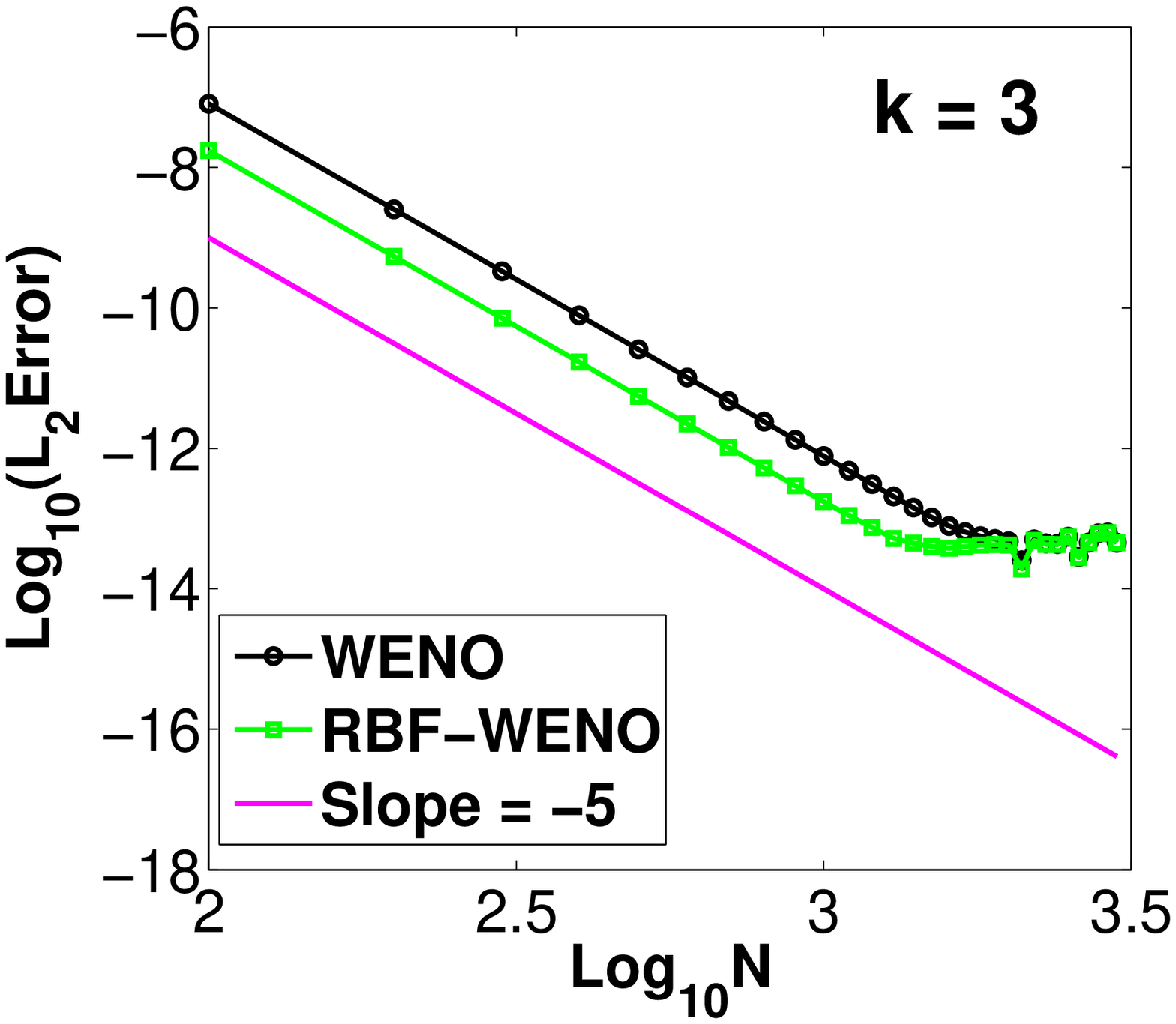}
\end{center}
\caption{(Color online).  $L_2$ errors versus $N$ on logarithmic scale with $k = 3$. Left: ENO method (blue) and RBF-ENO method (red). Right: WENO method (black) and RBF-WENO method (green).}
\label{test2}
\end{figure}

First we check the reconstruction error by the RBF-ENO and RBF-WENO method applied to a smooth function $u(x) = \sin(\pi  x)$ with $k = 2$ and $k=3$ to confirm the desired order of accuracy. Figures \ref{test1} and \ref{test2} show the $L_2$ errors versus $N$ on logarithmic scale. We see that the accuracy of the RBF-ENO is always one order higher than the accuracy of the regular ENO. The RBF-WENO method yields the same order of accuracy as the regular WENO method, 
but it is still better than the regular WENO in terms of accuracy because each reconstruction of the RBF-ENO is more accurate than the reconstruction by the regular ENO.

\subsubsection{Example 1}
We solve the smooth linear scalar equation given by the following advection equation for $x \in [-1,1]$
\begin{eqnarray}
u_t + u_x = 0, \quad  t > 0,
\label{example1}
\end{eqnarray}
with the initial condition $u (x,0) = u_0(x) = \sin(\pi x)$ and the periodic boundary condition. The CFL condition is given by $\Delta t \le C \Delta x$ with $C = 0.1$. 

Tables \ref{table0} and \ref{table1} show the $L_1, L_2$ and $L_\infty$ errors for each method at the final time $T = 0.5$ with $k = 2$ and $k = 3$ respectively. Consider the case for $k=2$. Since we test the convergence for a smooth problem, the switching condition is not applied. Let the WENO-JS denote the WENO method with the original smoothness indicators developed by Jiang and Shu \cite{WENO}. It is clear that for $k = 2$, the RBF-ENO method has almost $3$rd order convergence while the regular ENO method yields the convergence of $2$nd order or less. Also we observe that for $k = 2$ the RBF-ENO is even better than the WENO-JS in terms of accuracy while the rates of convergence are similar. The RBF-WENO with $k = 2$ is better than the WENO-JS in terms of accuracy and it achieves a higher order than the WENO-JS. Since both the RBF-ENO and the RBF-WENO have $3$rd order convergence in this case, it is interesting to see how similar they are due to the fact that we already take advantage of all the information at each step. We see this again in the following example.

\begin{table}[h]
\renewcommand{\arraystretch}{3.5}
\caption{$L_1, L_2$ and $L_\infty$ errors for the linear advection equation, \eqref{example1} with $k = 2$ at $T = 0.5$. }
\begin{center} \footnotesize
\renewcommand{\arraystretch}{1.1}
\begin{tabular}{|c|c|c|c|c|c|c|c|} 
\hline  
Method & N & $L_1$ error & $L_1$ order & $L_2$ error & $L_2$ order & $L_\infty$ error & $L_\infty$ order\\ 
\hline 
       & 10 & 1.09E-1 &   --   & 1.38E-1 &    --  & 2.18E-1 & --     \\  
       & 20 & 4.59E-2 & 1.2509 & 5.30E-2 & 1.3806 & 9.39E-2 & 1.2196 \\  
ENO    & 40 & 1.37E-2 & 1.7413 & 1.78E-2 & 1.5712 & 4.03E-2 & 1.2190 \\
k = 2  & 80 & 3.80E-3 & 1.8518 & 5.69E-3 & 1.6481 & 1.68E-2 & 1.2603 \\ 
       & 160& 1.02E-3 & 1.8999 & 1.80E-3 & 1.6597 & 6.91E-3 & 1.2848 \\ 
       & 320& 2.70E-4 & 1.9223 & 5.69E-4 & 1.6612 & 2.81E-3 & 1.2963 \\  
\hline
           & 10 & 1.76E-2 &    --  & 2.31E-2 &    --  & 4.17E-2 &    --  \\  
           & 20 & 2.47E-3 & 2.8362 & 2.64E-3 & 3.1294 & 3.61E-3 & 3.5316 \\  
RBF-ENO    & 40 & 3.17E-4 & 2.9628 & 3.43E-4 & 2.9443 & 4.78E-4 & 2.9163 \\
k = 2      & 80 & 4.05E-5 & 2.9693 & 4.42E-5 & 2.9581 & 6.25E-5 & 2.9424 \\ 
           & 160& 5.17E-6 & 2.9699 & 5.60E-6 & 2.9807 & 7.97E-6 & 2.9649 \\ 
           & 320& 6.51E-7 & 2.9896 & 7.05E-7 & 2.9891 & 1.00E-6 & 2.9894 \\
\hline
           & 10 & 8.94E-2 &   --   & 1.07E-1 &    --  & 1.69E-1 &   --   \\  
           & 20 & 2.90E-2 & 1.6239 & 3.23E-2 & 1.7370 & 5.47E-2 & 1.6324 \\  
WENO-JS    & 40 & 4.80E-3 & 2.5928 & 6.38E-3 & 2.3577 & 1.37E-2 & 1.9969 \\
k = 2      & 80 & 6.42E-4 & 2.9046 & 9.42E-4 & 2.7448 & 2.60E-3 & 2.3957 \\ 
           & 160& 7.79E-5 & 3.0433 & 1.26E-4 & 2.9420 & 3.96E-4 & 2.7176 \\ 
           & 320& 9.54E-6 & 3.0289 & 1.52E-5 & 3.0033 & 5.24E-5 & 2.9178 \\
\hline
           & 10 & 2.20E-2 &    --  & 2.27E-2 &    --  & 3.52E-2 &    --  \\  
           & 20 & 2.65E-3 & 3.0539 & 2.74E-3 & 3.0454 & 3.74E-3 & 3.2320 \\  
RBF-WENO   & 40 & 3.27E-4 & 3.0177 & 3.58E-4 & 2.9395 & 5.08E-4 & 2.8830 \\
k = 2      & 80 & 4.05E-5 & 3.0133 & 4.50E-5 & 2.9912 & 6.61E-5 & 2.9411 \\ 
           & 160& 5.09E-6 & 2.9932 & 5.63E-6 & 3.0001 & 8.27E-6 & 2.9998 \\ 
           & 320& 6.39E-7 & 2.9934 & 7.03E-7 & 3.0007 & 1.00E-6 & 3.0465 \\
\hline
\end{tabular}
\end{center} 
\label{table0}
\end{table}

\begin{table}[h]
\renewcommand{\arraystretch}{3.5}
\caption{$L_1, L_2$ and $L_\infty$ errors for the linear advection equation, \eqref{example1} with $k = 3$ at $T = 0.5$. }
\begin{center} \footnotesize
\renewcommand{\arraystretch}{1.1}
\begin{tabular}{|c|c|c|c|c|c|c|c|} 
\hline  
Method & N & $L_1$ error & $L_1$ order & $L_2$ error & $L_2$ order & $L_\infty$ error & $L_\infty$ order\\ 
\hline 
       & 10 & 2.88E-2 &   --   & 2.50E-2 &    --  & 3.59E-2 &   --   \\  
       & 20 & 2.78E-3 & 3.0390 & 3.03E-3 & 3.0547 & 4.45E-3 & 3.0135 \\  
ENO    & 40 & 3.36E-4 & 3.0475 & 3.68E-4 & 3.0316 & 5.47E-4 & 3.0226 \\
k = 3  & 80 & 4.12E-5 & 3.0271 & 4.54E-5 & 3.0177 & 6.76E-5 & 3.0177 \\ 
       & 160& 5.10E-6 & 3.0143 & 5.65E-6 & 3.0096 & 8.53E-6 & 2.9864 \\ 
       & 320& 6.34E-7 & 3.0073 & 7.03E-7 & 3.0050 & 1.06E-6 & 3.0091 \\  
\hline
           & 10 & 1.76E-2 &   --   & 1.93E-2 &    --  & 2.88E-2 &   --   \\  
           & 20 & 1.91E-3 & 3.1989 & 2.36E-3 & 3.0323 & 4.18E-3 & 2.7936 \\  
RBF-ENO    & 40 & 1.44E-4 & 3.7303 & 2.08E-4 & 3.5015 & 4.98E-4 & 3.0674 \\
k = 3      & 80 & 8.79E-6 & 4.0371 & 1.56E-5 & 3.7343 & 5.09E-5 & 3.2902 \\ 
           & 160& 4.81E-7 & 4.1923 & 1.05E-6 & 3.8872 & 4.59E-6 & 3.4715 \\ 
           & 320& 2.76E-8 & 4.1204 & 7.27E-8 & 3.8628 & 4.23E-7 & 3.4387 \\  
\hline
           & 10 & 9.57E-3 &   --   & 1.12E-2 &    --  & 1.60E-2 &   --   \\  
           & 20 & 3.99E-4 & 4.5829 & 4.62E-4 & 4.5948 & 7.80E-4 & 4.3626 \\  
WENO-JS    & 40 & 1.18E-5 & 5.0807 & 1.38E-5 & 5.0646 & 2.47E-5 & 4.9824 \\
k = 3      & 80 & 3.70E-7 & 4.9946 & 4.28E-7 & 5.0111 & 7.82E-7 & 4.9802 \\ 
           & 160& 1.34E-8 & 4.7841 & 1.51E-8 & 4.8232 & 2.67E-8 & 4.8696 \\ 
           & 320& 6.56E-10& 4.3553 & 7.30E-10& 4.3717 & 1.13E-9 & 4.5696 \\ 
\hline      
		   & 10 & 2.69E-3 &   --   & 2.93E-3 &    --  & 4.19E-3 &   --   \\  
           & 20 & 8.92E-5 & 4.9139 & 1.05E-4 & 4.8006 & 1.94E-4 & 4.4291 \\  
RBF-WENO   & 40 & 2.53E-6 & 5.1401 & 3.00E-6 & 5.1268 & 6.26E-6 & 4.9568 \\
k = 3      & 80 & 7.52E-8 & 5.0710 & 8.56E-8 & 5.1353 & 1.55E-7 & 5.3351 \\ 
           & 160& 2.35E-9 & 5.0009 & 2.63E-9 & 5.0245 & 4.95E-9 & 4.9693 \\ 
           & 320& 7.39E-11& 4.9902 & 8.32E-11& 4.9819 & 1.76E-10& 4.8097 \\  
\hline
\end{tabular}
\end{center} 
\label{table1}
\end{table}

\subsubsection{Example 2}
We consider the same advection equation \eqref{example1} but with the discontinuous initial condition 
\begin{eqnarray}
   u(x,0) = -\mbox{sgn}(x),
   \label{example2} 
\end{eqnarray}
and the boundary condition $u(-1,t) = 1, t>0$. With this example, we check how the RBF-ENO solution behaves near the discontinuity. 

\begin{figure}[h]
\begin{center}
\includegraphics[width=0.48\textwidth]{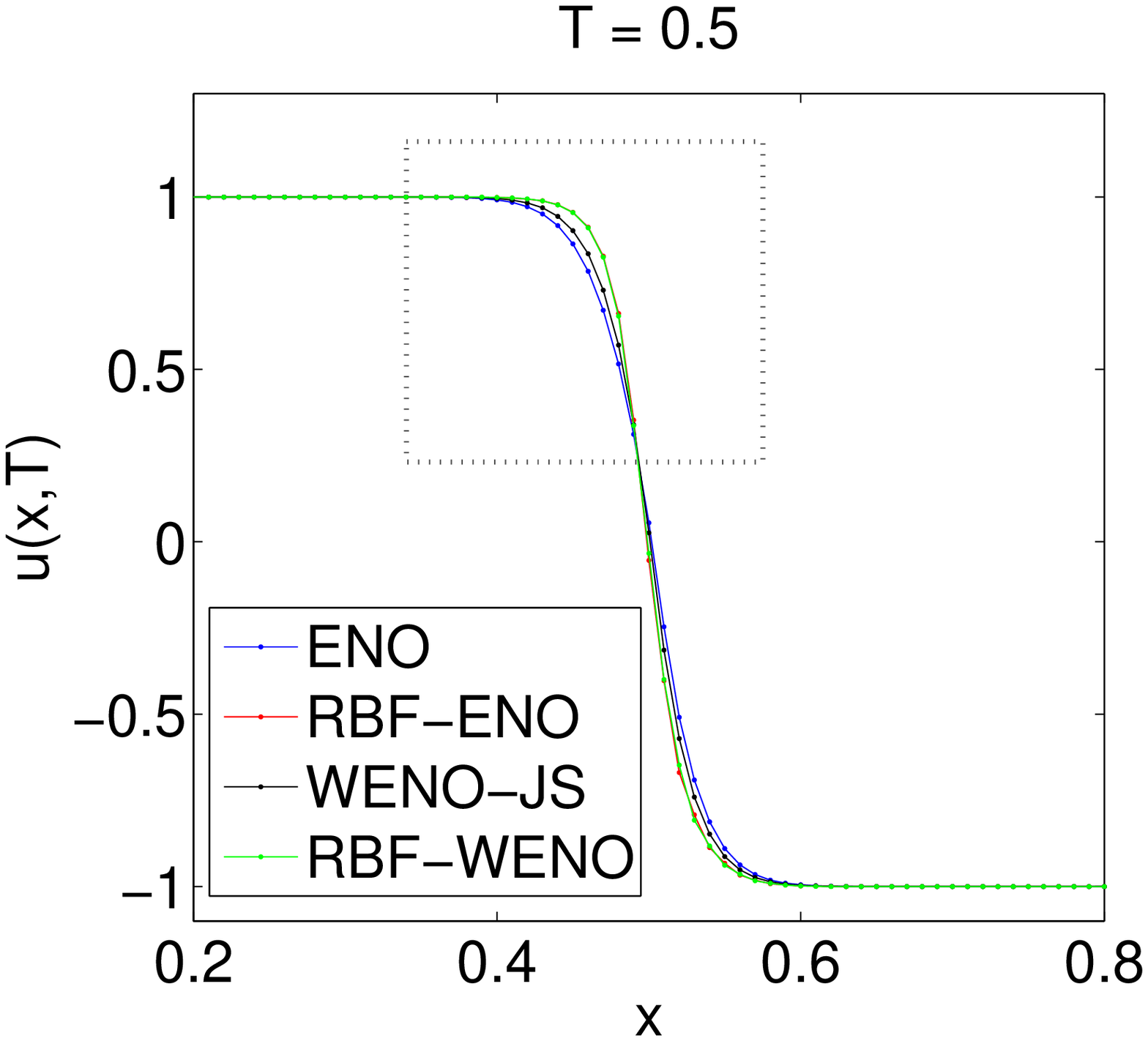}
\includegraphics[width=0.48\textwidth]{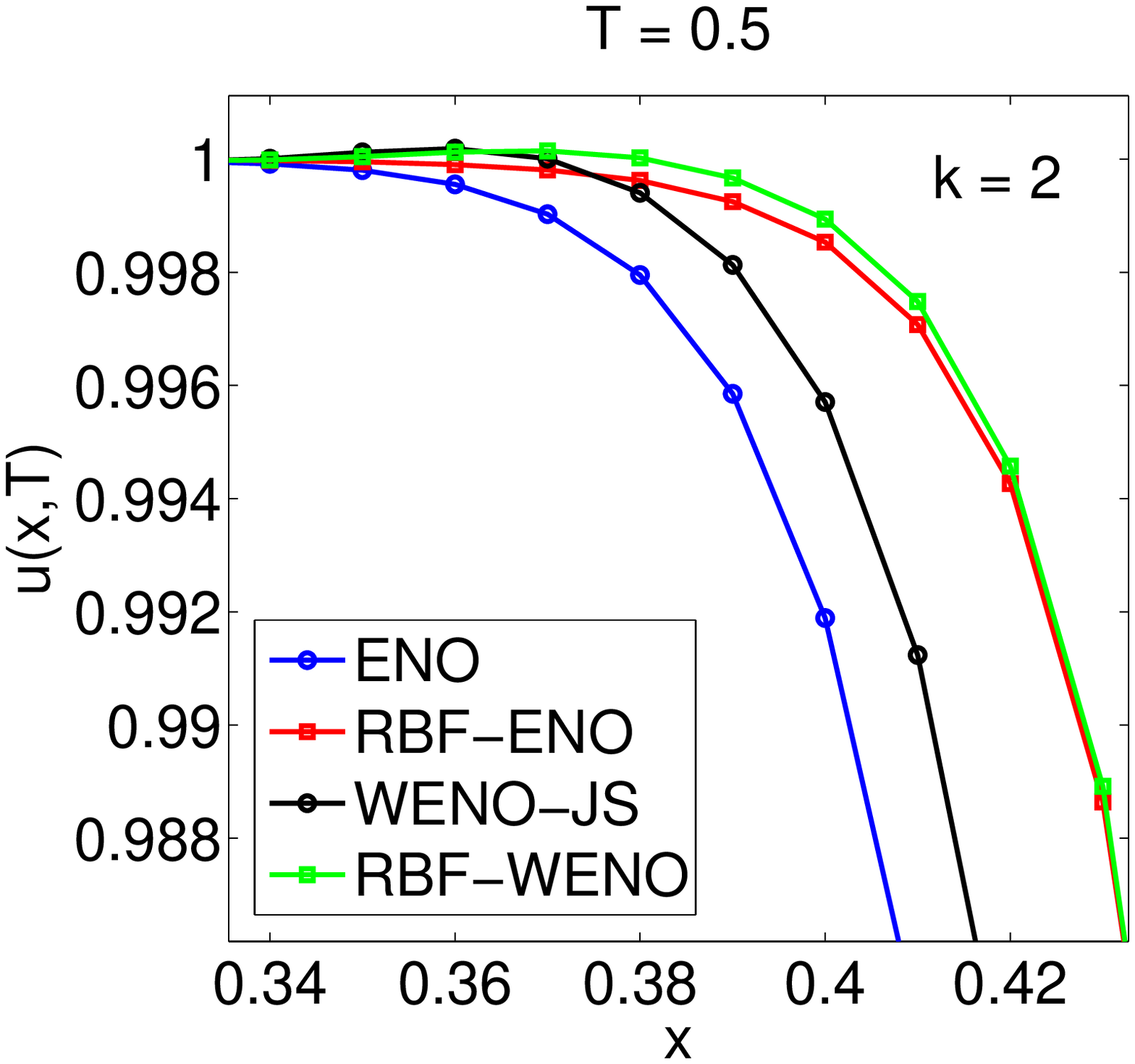}
\includegraphics[width=0.48\textwidth]{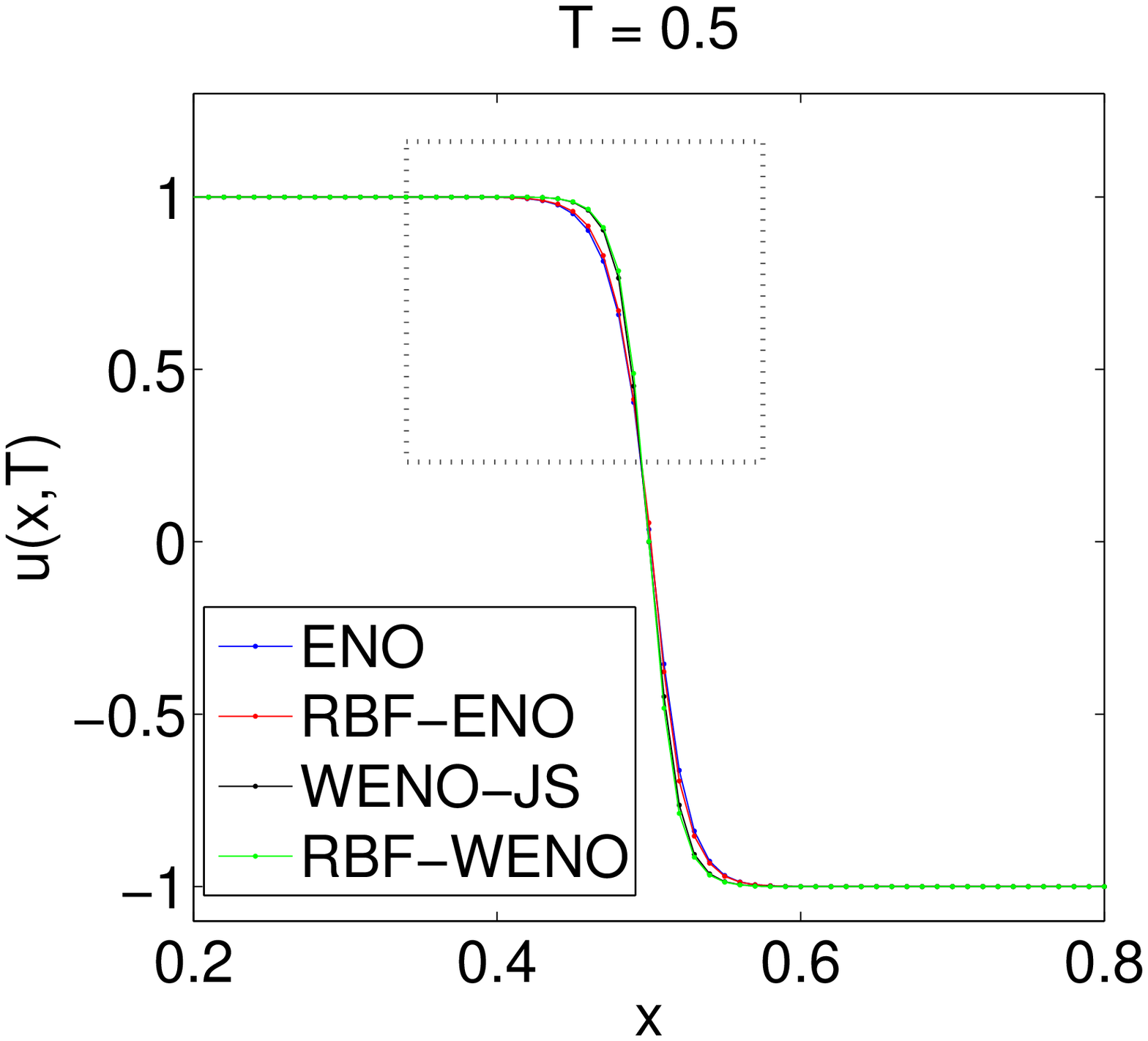}
\includegraphics[width=0.48\textwidth]{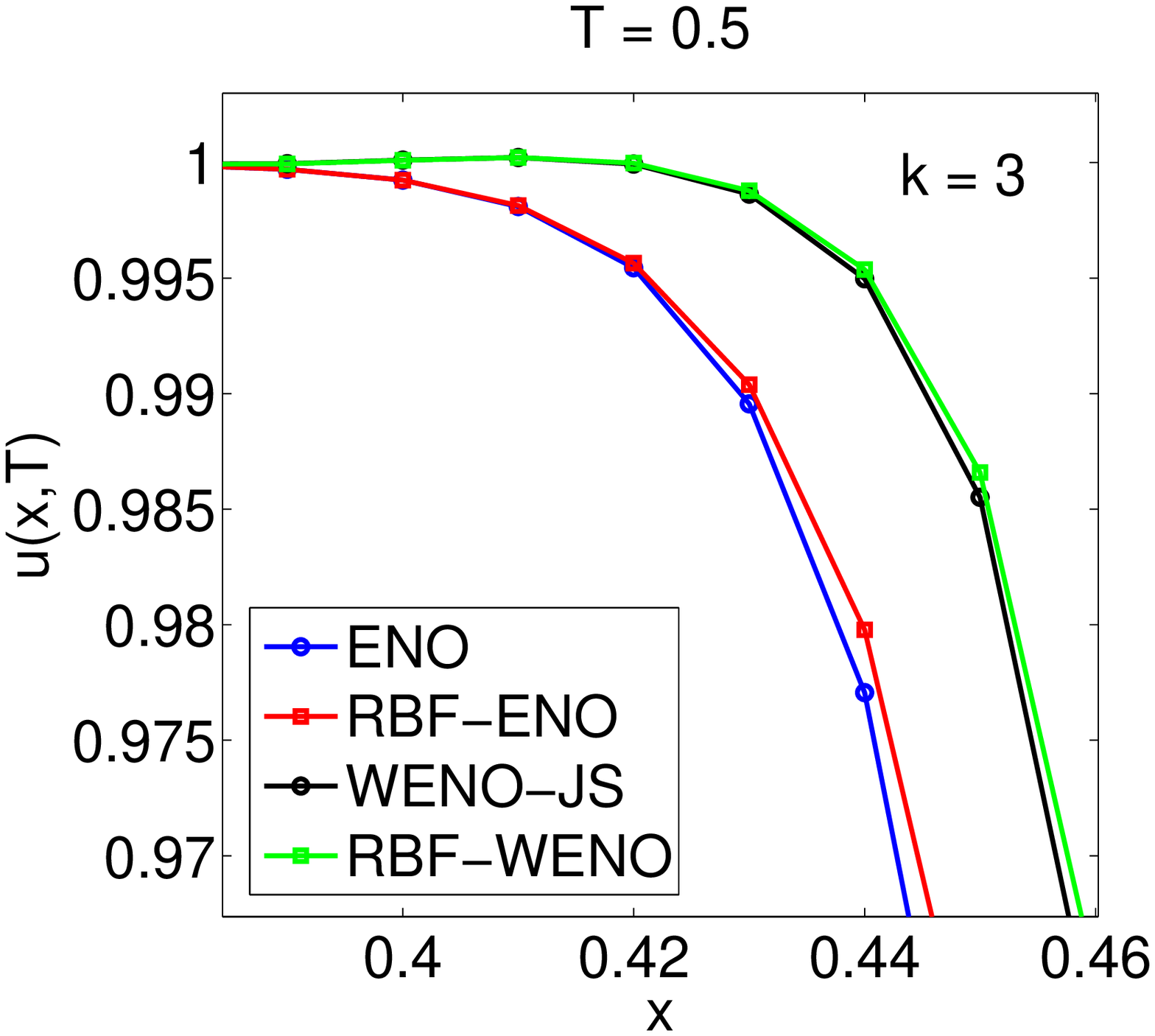}
\end{center}
\caption{(Color online). Solutions to \eqref{example1} at $T = 0.5$ with the discontinuous initial condition, \eqref{example2} for the ENO (blue), RBF-ENO (red), WENO-JS (black) and RBF-WENO (green) methods with $k = 2$ (top) and $k = 3$ (bottom). $N = 200$.}
\label{figure2}
\end{figure}
Figure \ref{figure2} shows the solution profiles at $T = 0.5$ by each method with $N = 200$. The top two figures show the solutions with $k = 2$ and the bottom two figures with $k = 3$. As shown in the figures, the RBF-ENO solutions for both $k = 2$ and $k = 3$ are non-oscillatory. For the case of $k = 2$, the RBF-ENO solution is superior to the regular ENO and WENO-JS solutions. The ENO or RBF-ENO methods give non-oscillatory solutions due to the adaptive stencil at each reconstruction step while the WENO-JS or RBF-WENO methods have no such guarantee since they use a combination of all the possible ENO reconstruction including the oscillatory ones. This is why when the RBF-ENO, WENO-JS and RBF-WENO methods have the similar $3$rd order convergence for $k=2$ and the RBF-ENO gives the best non-oscillatory profile.

For the case of $k=3$, the RBF-ENO is better than the regular ENO because it has convergence of one order higher, while the RBF-WENO is slightly better than the WENO-JS in terms of accuracy and both of them have $5$th order convergence.

\subsubsection{Example 3}
\begin{figure}[h]
\begin{center}
\includegraphics[width=0.496\textwidth]{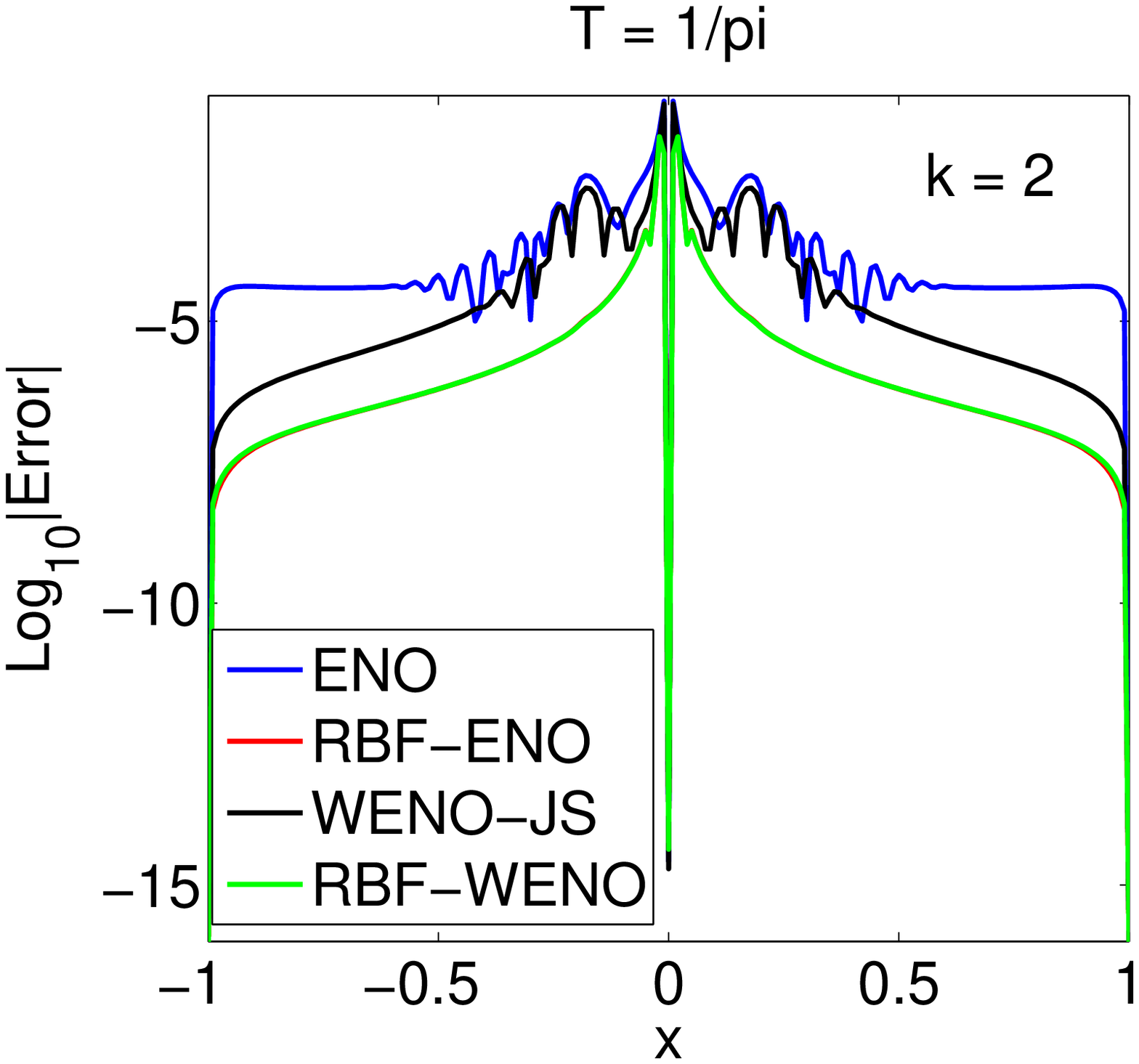}
\includegraphics[width=0.496\textwidth]{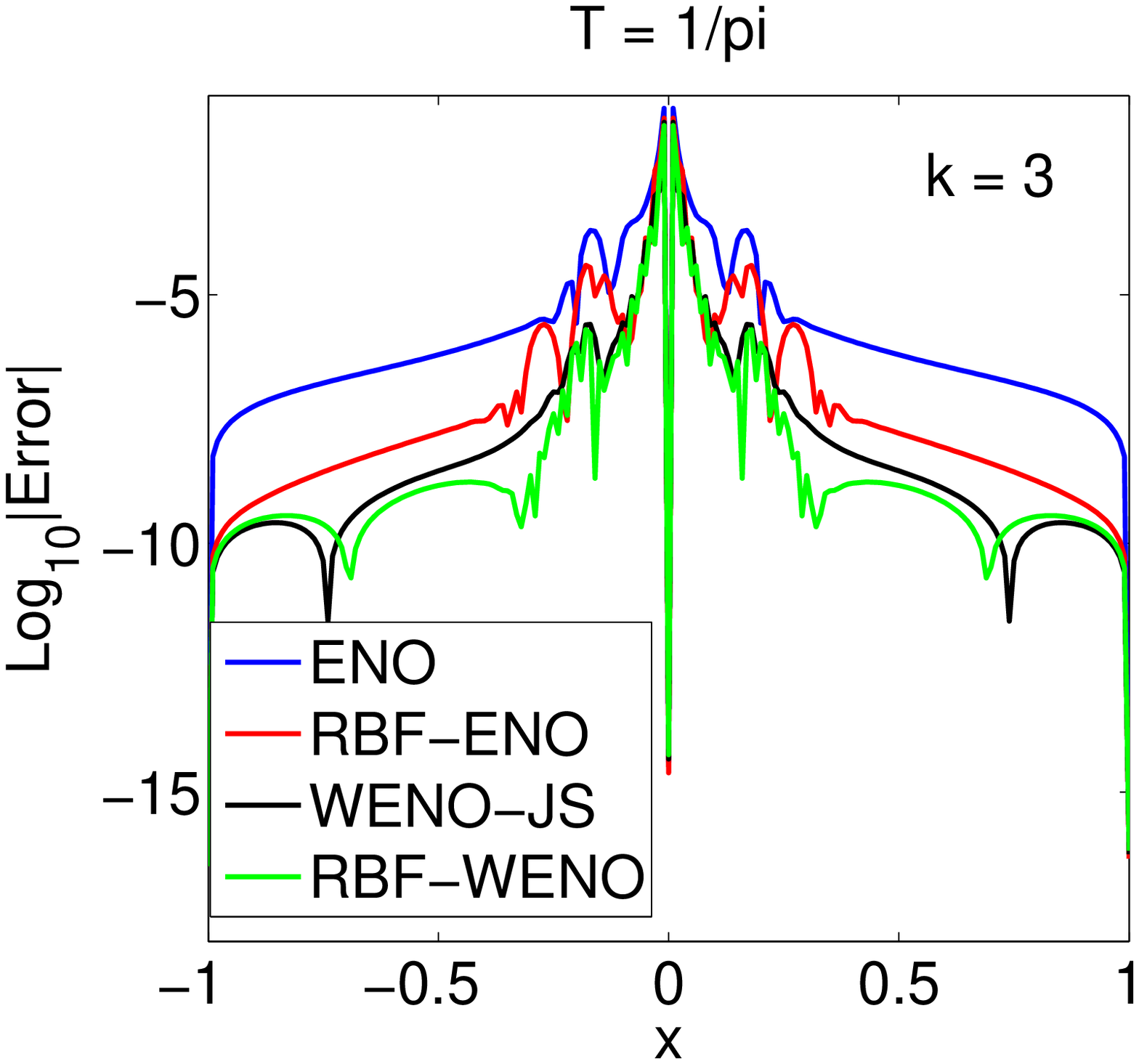}
\end{center}
\caption{(Color online). Pointwise errors for Burgers equation at $T = 1/\pi$ by the ENO (blue), RBF-ENO (red), WENO-JS (black) and RBF-WENO (green) methods.  $N = 200$. Left: $k = 2$. Right: $k = 3$. }
\label{figure3}
\end{figure}

We consider the Burgers' equation for $x \in [-1,1]$ 
\begin{eqnarray}
   \label{example3}
   u_t + \left(\frac{1}{2} u^2 \right)_x &=& 0,  \quad t >0, 
  \\
   \label{example3_example}
   u(x,0) &=& -\sin(\pi x). 
\end{eqnarray}
%
For this problem, a shock forms at $(x,t) = (0, {1/\pi})$. 
We first check if each method achieves its desired order when the solution is still smooth at $T= 0.2$. Tables \ref{table2} and \ref{table3} show the $L_1, L_2$ and $L_\infty$ errors by the ENO, RBF-ENO, WENO-JS and RBF-WENO methods for $k = 2$ and $k = 3$. Figure \ref{figure3} shows the pointwise errors of the solution at $T = {1/ \pi}$ with $k = 2$ and $k = 3$. The RBF-ENO method and RBF-WENO method yield almost the same accurate results due to the fact that we already took advantage of all the available information. For the case of $k = 3$, the RBF-WENO method yields the most accurate results. From this figure, we observe that the RBF-ENO method gives smaller pointwise errors than the regular ENO method, which confirms that the RBF-ENO method much improves the original ENO accuracy in the smooth area.

\begin{table}[h]
\renewcommand{\arraystretch}{3.5}
\caption{$L_1, L_2$ and $L_\infty$ errors for the Burgers' equation, \eqref{example3} with $k = 2$ at $T = 0.2$. }
\begin{center} \footnotesize
\renewcommand{\arraystretch}{1.1}
\begin{tabular}{|c|c|c|c|c|c|c|c|} 
\hline  
Method & N & $L_1$ error & $L_1$ order & $L_2$ error & $L_2$ order & $L_\infty$ error & $L_\infty$ order \\
\hline
           & 10 & 9.58E-2 &   --   & 1.29E-1 &    --  & 2.84E-1 &    --  \\  
           & 20 & 3.01E-2 & 1.6687 & 4.54E-2 & 1.5086 & 1.12E-1 & 1.3319 \\  
ENO        & 40 & 9.03E-3 & 1.7392 & 1.49E-2 & 1.6078 & 5.68E-2 & 0.9903 \\
k = 2      & 80 & 2.49E-3 & 1.8552 & 4.33E-3 & 1.7826 & 1.64E-2 & 1.7885 \\ 
           & 160& 6.77E-4 & 1.8805 & 1.23E-3 & 1.8095 & 4.17E-3 & 1.8017 \\ 
           & 320& 1.78E-4 & 1.9243 & 3.51E-4 & 1.8137 & 1.52E-3 & 1.6295 \\  
\hline
           & 10 & 5.40E-2 &   --   & 9.82E-2 &    --  & 2.18E-1 &    --  \\  
           & 20 & 9.56E-3 & 2.4993 & 2.15E-2 & 2.1886 & 7.47E-2 & 1.5464 \\  
RBF-ENO    & 40 & 1.46E-3 & 2.7042 & 3.55E-3 & 2.6005 & 1.65E-2 & 2.1781 \\
k = 2      & 80 & 1.85E-4 & 2.9794 & 4.96E-4 & 2.8390 & 2.58E-3 & 2.6762 \\ 
           & 160& 2.28E-5 & 3.0257 & 6.20E-5 & 3.0005 & 3.69E-4 & 2.8071 \\ 
           & 320& 2.78E-6 & 3.0373 & 7.55E-6 & 3.0378 & 4.51E-5 & 3.0299 \\
\hline
           & 10 & 7.45E-2 &   --   & 1.13E-1 &    --  & 2.83E-1 &    --   \\  
           & 20 & 2.24E-2 & 1.7337 & 3.77E-2 & 1.5833 & 1.02E-1 & 1.4739  \\  
WENO-JS    & 40 & 4.46E-3 & 2.3269 & 8.17E-3 & 2.2081 & 2.99E-2 & 1.7701  \\
k = 2      & 80 & 6.31E-4 & 2.8234 & 1.13E-3 & 2.8507 & 4.12E-3 & 2.8594  \\ 
           & 160& 8.07E-5 & 2.9668 & 1.44E-4 & 2.9748 & 5.13E-4 & 3.0064  \\ 
           & 320& 1.01E-5 & 2.9943 & 1.79E-5 & 3.0032 & 6.28E-5 & 3.0322  \\
\hline 
           & 10 & 5.44E-2 &   --   & 9.82E-2 &    --  & 2.18E-1 &    --  \\  
           & 20 & 9.67E-3 & 2.4938 & 2.15E-2 & 2.1875 & 7.55E-2 & 1.5317 \\  
RBF-WENO   & 40 & 1.44E-3 & 2.7471 & 3.54E-3 & 2.6067 & 1.66E-2 & 2.1827 \\
k = 2      & 80 & 1.86E-4 & 2.9504 & 4.95E-4 & 2.8386 & 2.60E-3 & 2.6757 \\ 
           & 160& 2.27E-5 & 3.0346 & 6.18E-5 & 3.0009 & 3.66E-4 & 2.8285 \\ 
           & 320& 2.77E-6 & 3.0356 & 7.53E-6 & 3.0364 & 4.50E-5 & 3.0243 \\
\hline 
\end{tabular}
\end{center} 
\label{table2}
\end{table}

\begin{table}[h]
\renewcommand{\arraystretch}{3.5}
\caption{$L_1, L_2$ and $L_\infty$ errors for the Burgers' equation, \eqref{example3} with $k = 3$ at $T = 0.2$. }
\begin{center} \footnotesize
\renewcommand{\arraystretch}{1.1}
\begin{tabular}{|c|c|c|c|c|c|c|c|} 
\hline
Method & N & $L_1$ error & $L_1$ order & $L_2$ error & $L_2$ order & $L_\infty$ error & $L_\infty$ order \\
\hline
           & 10 & 4.32E-2 &   --   & 8.33E-2 &    --  & 2.49E-1 &  --     \\  
           & 20 & 9.68E-3 & 2.1580 & 2.02E-2 & 2.0445 & 7.24E-2 & 1.7852  \\  
ENO        & 40 & 1.47E-3 & 2.7194 & 3.11E-3 & 2.6989 & 1.44E-2 & 2.3289  \\
k = 3      & 80 & 2.23E-4 & 2.7201 & 4.99E-4 & 2.6387 & 2.47E-3 & 2.5402  \\ 
           & 160& 3.11E-5 & 2.8450 & 7.05E-5 & 2.8253 & 3.63E-4 & 2.7703  \\ 
           & 320& 4.31E-6 & 2.8485 & 9.60E-6 & 2.8754 & 4.45E-5 & 3.0271  \\  
\hline         
           & 10 & 3.45E-2 &   --   & 6.61E-2 &    --  & 1.94E-1 &   --   \\  
           & 20 & 7.76E-3 & 2.1550 & 1.91E-2 & 1.7864 & 6.57E-2 & 1.5676 \\  
RBF-ENO    & 40 & 1.24E-3 & 2.6451 & 4.47E-3 & 2.0061 & 2.90E-2 & 1.1769 \\
k = 3      & 80 & 8.64E-5 & 3.8429 & 3.71E-4 & 3.6849 & 3.09E-3 & 3.2325 \\ 
           & 160& 8.39E-6 & 3.3645 & 3.37E-5 & 3.4584 & 3.83E-4 & 3.0129 \\ 
           & 320& 6.14E-7 & 3.7730 & 1.66E-6 & 4.3410 & 1.75E-5 & 4.4509 \\ 
\hline
           & 10 & 3.25E-2 &   --   & 7.21E-2 &    --  & 2.21E-1 &   --    \\  
           & 20 & 4.24E-3 & 2.9386 & 1.33E-2 & 2.4336 & 5.83E-2 & 1.9182  \\  
WENO-JS    & 40 & 4.19E-4 & 3.3372 & 1.46E-3 & 3.1896 & 8.75E-3 & 2.7373  \\
k = 3      & 80 & 2.45E-5 & 4.0994 & 9.00E-5 & 4.0221 & 6.06E-4 & 3.8522  \\ 
           & 160& 9.42E-7 & 4.6994 & 3.39E-6 & 4.7292 & 2.50E-5 & 4.5983  \\ 
           & 320& 2.94E-8 & 5.0026 & 1.08E-7 & 4.9776 & 8.18E-7 & 4.9346  \\
\hline  
           & 10 & 3.45E-2 &   --   & 5.90E-2 &    --  & 1.55E-1 &   --   \\  
           & 20 & 3.77E-3 & 3.1931 & 9.28E-3 & 2.6681 & 3.87E-2 & 2.0069 \\  
RBF-WENO   & 40 & 3.17E-4 & 3.5739 & 9.66E-4 & 3.2642 & 5.06E-3 & 2.7900 \\
k = 3      & 80 & 1.86E-5 & 4.0856 & 5.60E-5 & 4.1075 & 3.66E-4 & 3.9367 \\ 
           & 160& 9.62E-7 & 4.2793 & 2.57E-6 & 4.4464 & 1.39E-5 & 4.7109 \\ 
           & 320& 2.86E-8 & 5.0698 & 8.21E-8 & 4.9685 & 4.67E-7 & 4.9024 \\ 
\hline 
\end{tabular}
\end{center} 
\label{table3}
\end{table}


\newpage
\subsubsection{Example 4 - system} 
\begin{figure}[h]
\begin{center}
\includegraphics[width=0.48\textwidth]{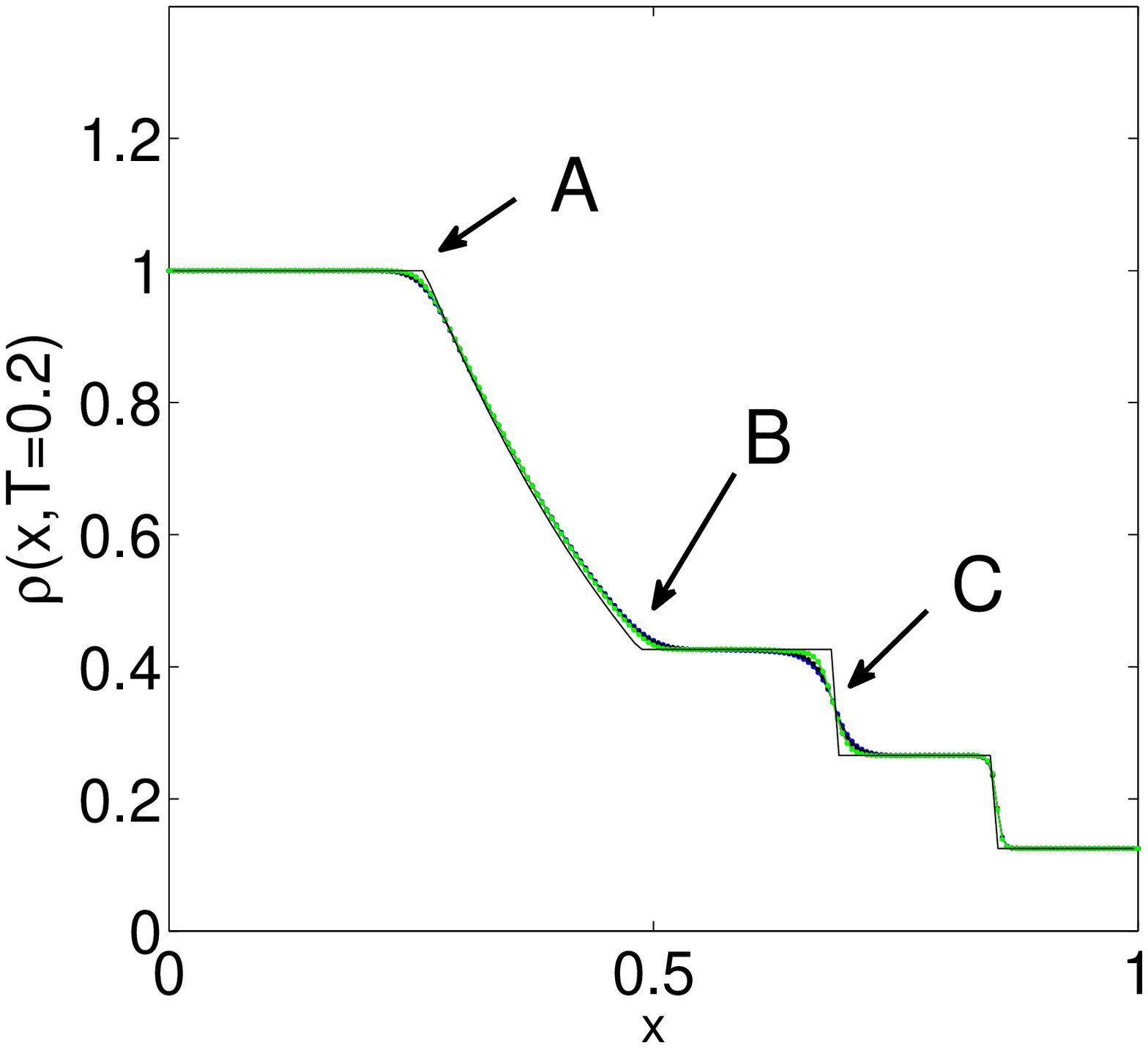}
\includegraphics[width=0.48\textwidth]{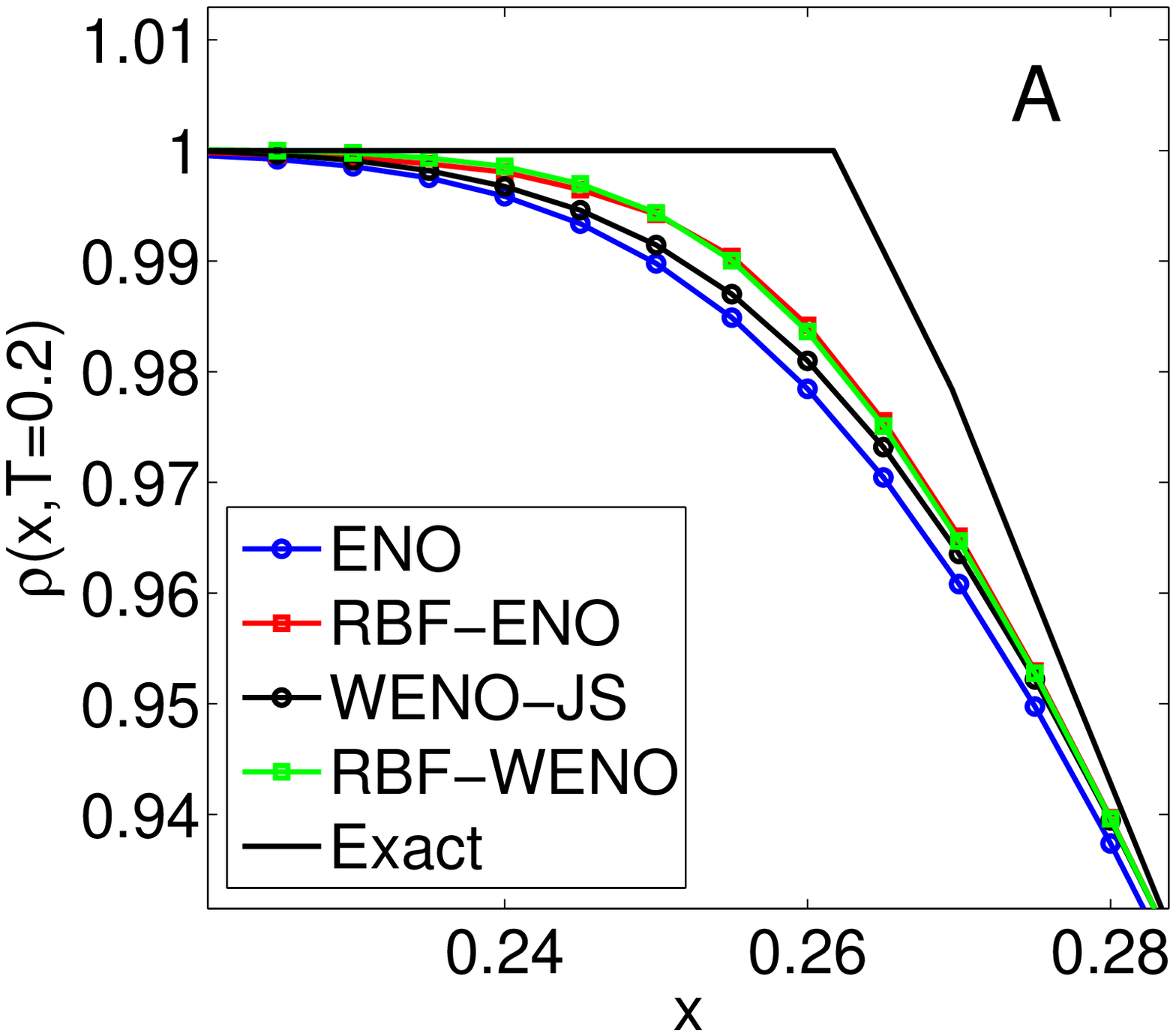}
\includegraphics[width=0.48\textwidth]{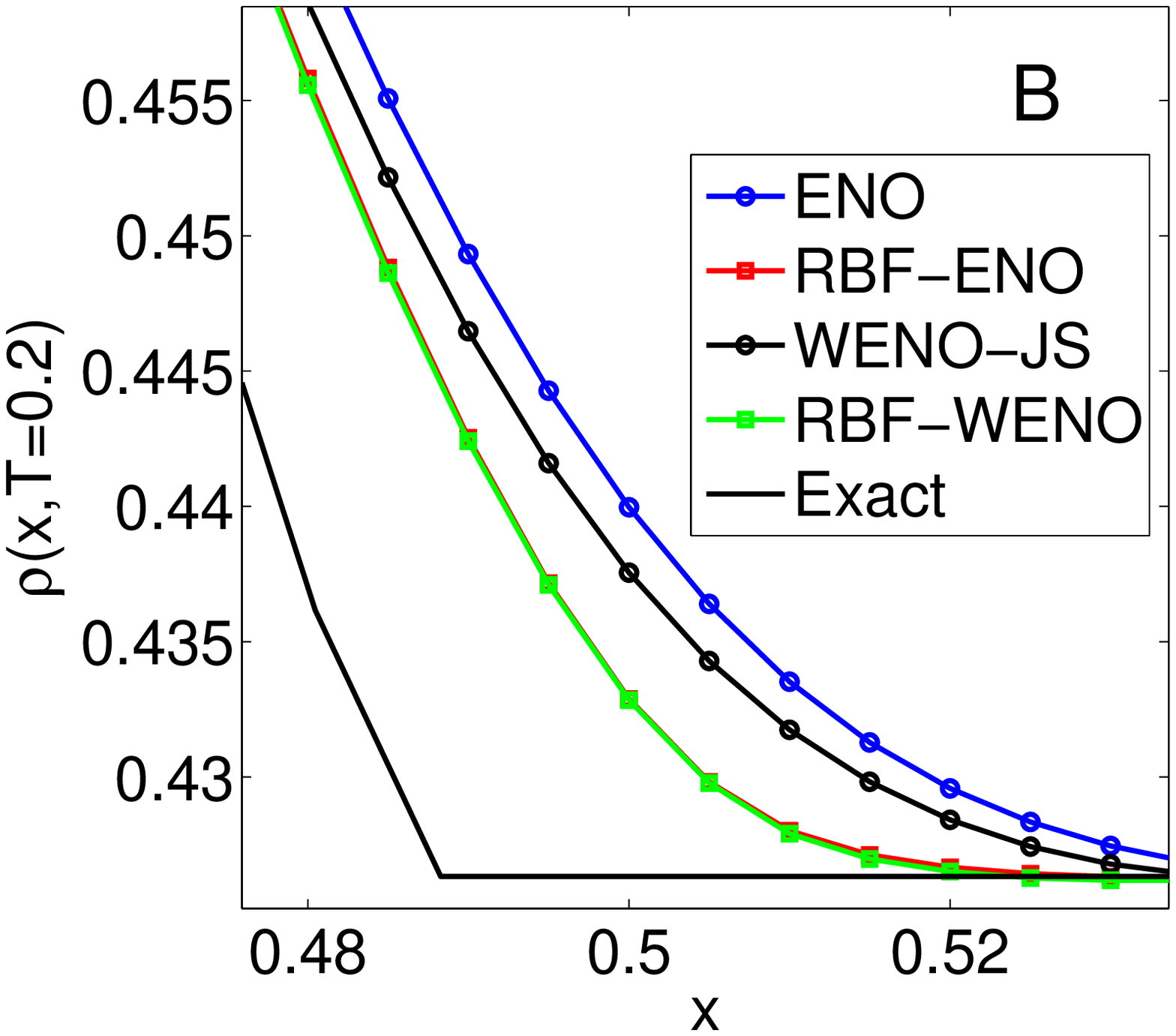}
\includegraphics[width=0.48\textwidth]{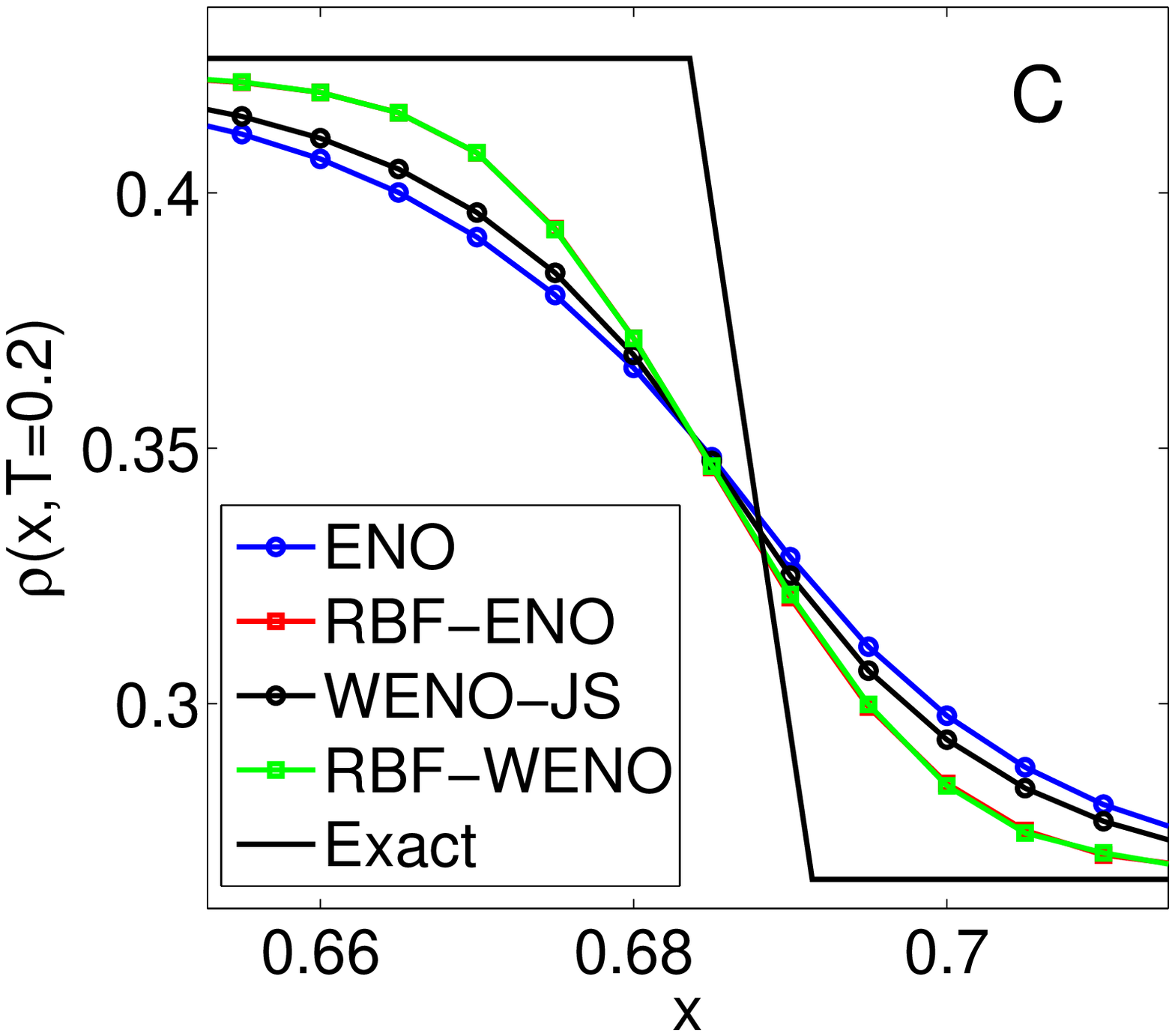}
\end{center}
\caption{(Color online). Density profile for Sod problem for the ENO (blue), RBF-ENO (red), WENO-JS (black), RBF-WENO (green) and the exact solution at $T = 0.2$ with $N = 600$.  $k = 2$.}
\label{sod}
\end{figure}

Finally we consider the one-dimensional Euler equations for gas dynamics
\begin{eqnarray}
    U_t + F(U)_x = 0, 
\label{euler}
\end{eqnarray}
where the conservative state vector $U$ and the flux function $F$ are given by 
$$
  U = (\rho, \rho u, E)^T, \quad F(U) = (\rho u, \rho u^2 + P, (E+P)u)^T.
$$
Here $\rho, u, P$ and $E$ denote density, velocity, pressure and total energy, respectively. The equation of state is given by $$P = (\gamma - 1) \left( E - {1\over 2} \rho u^2 \right), $$ where $\gamma = 1.4$ for the ideal gas. 
We consider the Sod shock tube problem with the initial condition 
$$
   (\rho, u, P ) = \left\{ \begin{array}{ll} (1, 0, 1) & x \le 0 \\ (0.125, 0, 0.1) & x > 0 \end{array} \right. .
$$
The Sod problem is solved with $N = 600$ and the CFL number $C = 0.1$. Figure \ref{sod} shows the density profile at $T = 0.2$ for $k = 2$ by each method. The exact solution is also provided. The exact solution is plotted on the $2000$ points. The top left figure shows the overall density profile and the other figures show the detailed density profiles for the regions of $A, B$ and $C$ indicated in the top left figure. In all areas of $A$, $B$ and $C$, the RBF-ENO method yields better (sharper) density profiles than the ENO and WENO methods, while it is almost the same as the RBF-WENO. We observe that the WENO solution is slightly oscillatory while the RBF-ENO and ENO solutions are not oscillatory. The detailed profiles near the shock are not given in the figure, but the density profile by the RBF-ENO method is slightly better than the regular ENO and WENO methods near the shock area as well. 

\begin{figure}[h]
\begin{center}
\includegraphics[width=0.48\textwidth]{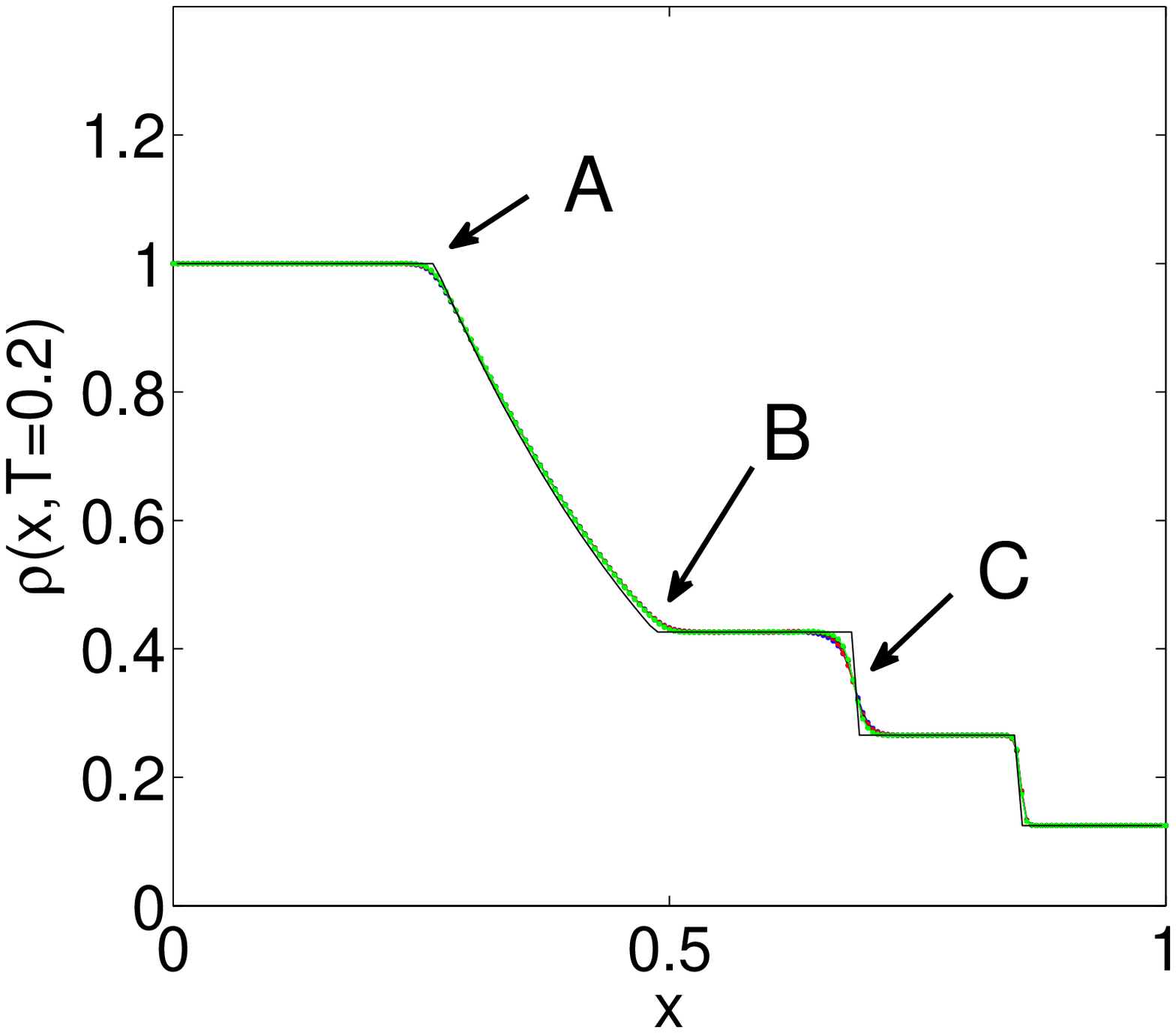}
\includegraphics[width=0.48\textwidth]{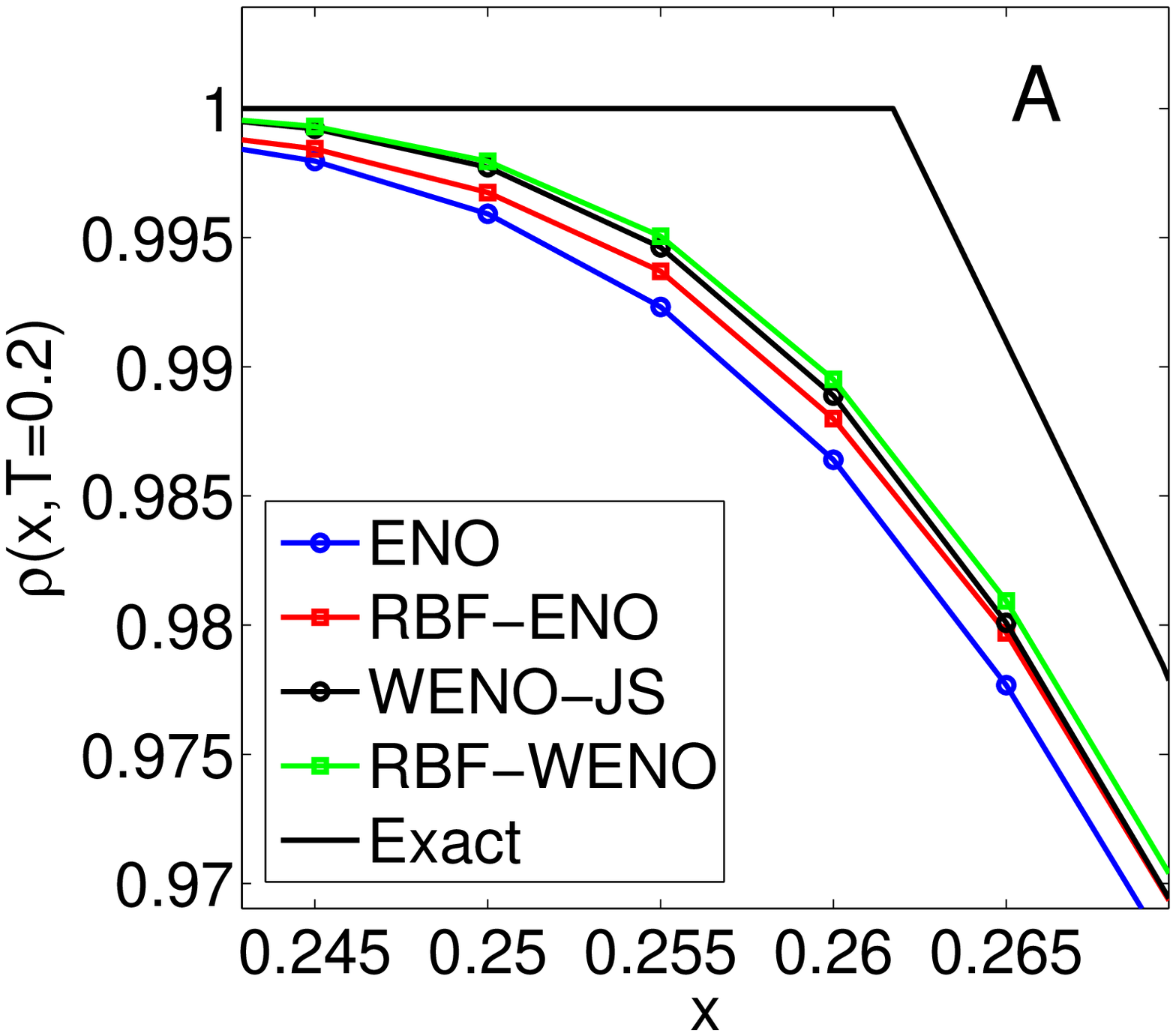}
\includegraphics[width=0.48\textwidth]{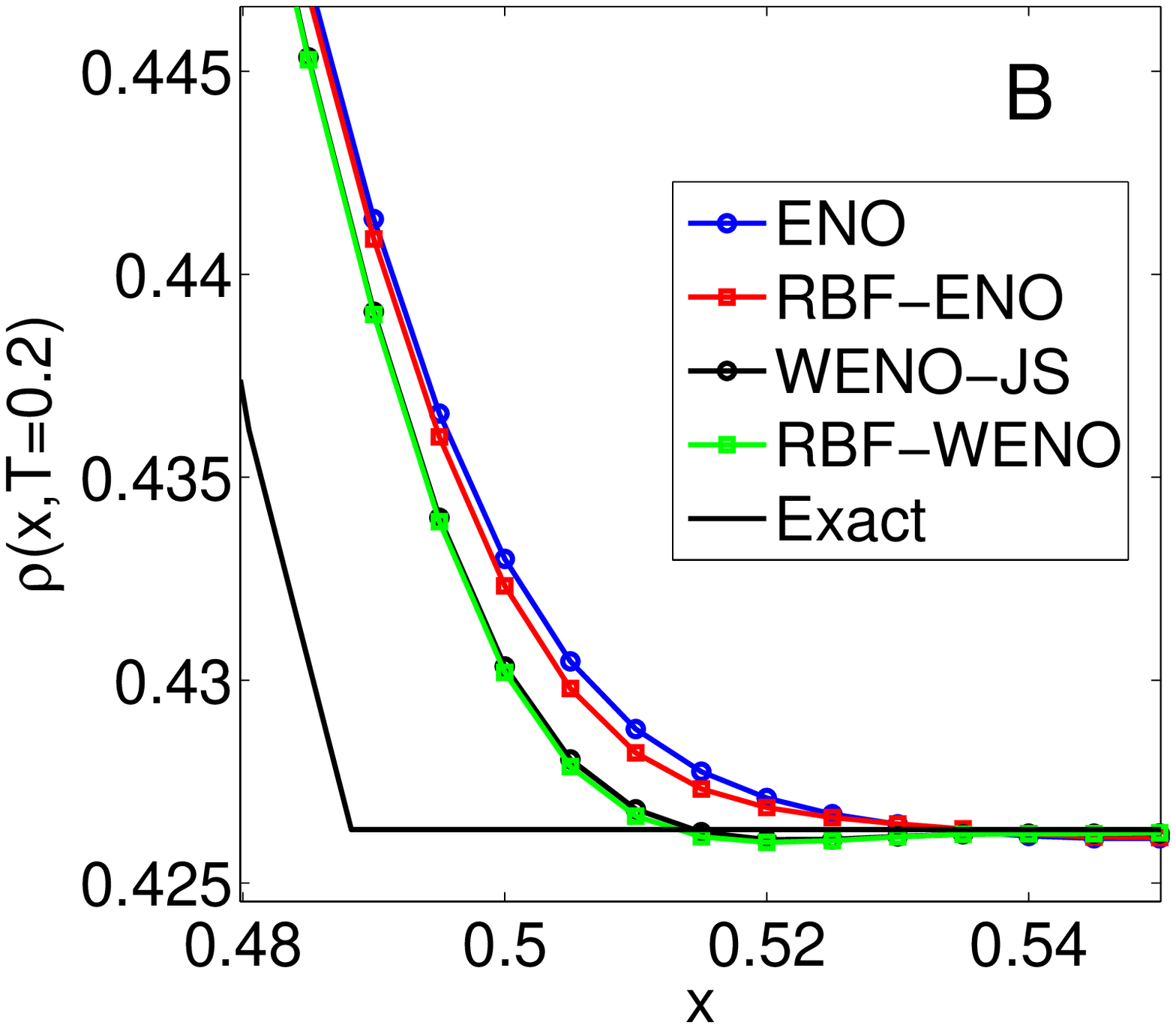}
\includegraphics[width=0.48\textwidth]{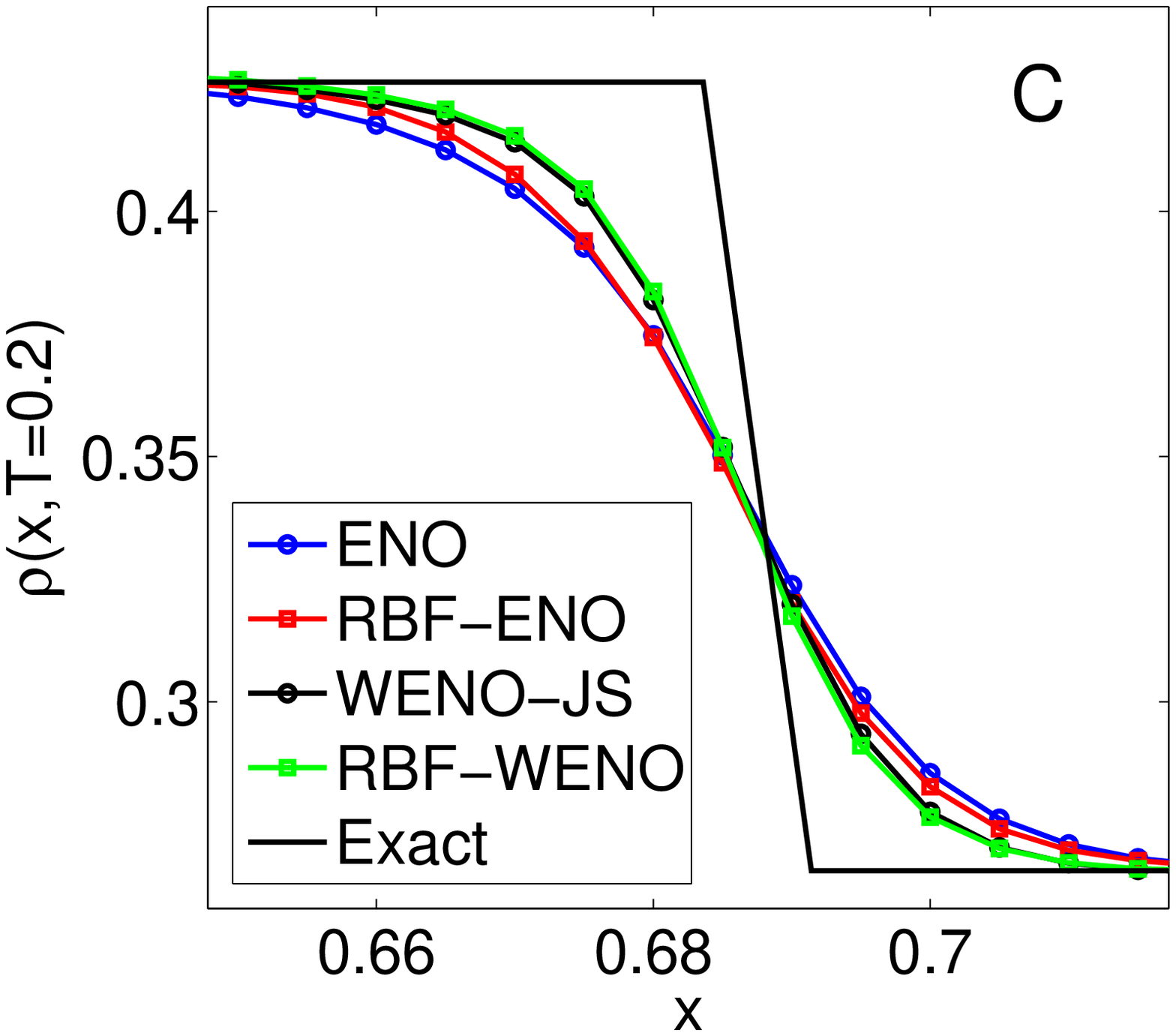}
\end{center}
\caption{(Color online). Density profile for Sod problem for the ENO (blue), RBF-ENO (red), WENO-JS (black), RBF-WENO (green) and the exact solution at $T = 0.2$ with $N = 600$.   $k = 3$.}
\label{sod2}
\end{figure}

Figure \ref{sod2} shows the same density profiles for $k = 3$. For this case, the RBF-ENO method yields $4$th order accuracy in the smooth area while the ENO and WENO-JS methods yield $3$rd order and $5$th order accuracy, respectively. Thus unlike the case of $k = 2$, it is reasonable to observe that the RBF-WENO method yields the best result among those three methods. The RBF-ENO solution is, however, still better than the regular ENO solution. Near the shock area, all three methods yield almost similar profiles but the RBF-WENO solution is sharper than the others and the RBF-ENO solution is sharper than the regular ENO method. 

\begin{center}
\small{
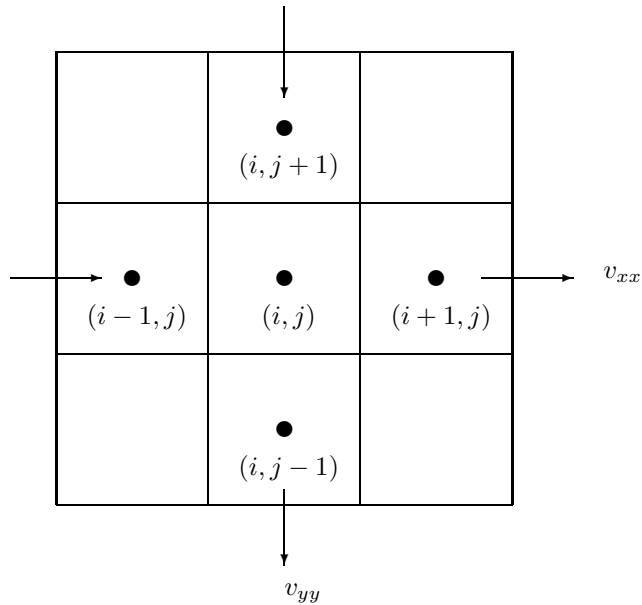
\begin{figure}
\setlength{\unitlength}{2cm}
\begin{picture}(6,4)(-3,-2)

\put(-1.5,0.5){\line(1,0){3}}
\put(-1.5,-0.5){\line(1,0){3}}
\put(-1.5,1.5){\line(1,0){3}}
\put(-1.5,-1.5){\line(1,0){3}}

\put(0.5,-1.5){\line(0,1){3}}
\put(-0.5,-1.5){\line(0,1){3}}
\put(1.5,-1.5){\line(0,1){3}}
\put(-1.5,-1.5){\line(0,1){3}}


\put(0,0){\circle*{0.1}}
	\put(-0.15,-0.3){$(i,j)$}
\put(0,1){\circle*{0.1}}
	\put(-0.3,0.7){$(i,j+1)$}
\put(0,-1){\circle*{0.1}}
	\put(-0.3,-1.3){$(i,j-1)$}
\put(1,0){\circle*{0.1}}
	\put(-1.3,-0.3){$(i-1,j)$}	
\put(-1,0){\circle*{0.1}}
	\put(0.7,-0.3){$(i+1,j)$}	
	
\put(1.3,0){\vector(1,0){0.6}}
\put(-1.8,0){\vector(1,0){0.6}}
\put(2.1,0){$v_{xx}$}

\put(0,1.8){\vector(0,-1){0.6}}
\put(0,-1.4){\vector(0,-1){0.5}}
\put(0,-2.1){$v_{yy}$}


%



%
%


%
%


\end{picture}
\caption{Five cell averages involved for the reconstruction at the cell boundaries  of $I_{i,j}$ for $k = 2$. }
\label{figure2d_grid}
\end{figure}
}
\end{center}

\break
\section{2D interpolation} 
To illustrate the non-polynomial reconstruction in 2D, we consider the case with $k=2$, for which at least ${1\over 2} k(k+1)=3$ cells are needed for the ENO reconstruction among the available $5$ cell averages, ${\bar v}_{i-1,j}$, ${\bar v}_{i,j}$, ${\bar v}_{i+1,j}$, ${\bar v}_{i,j-1}$ and ${\bar v}_{i,j+1}$ (see Figure \ref{figure2d_grid}). For simplicity, we use the uniform grid, i.e. $ \Delta x_i = \Delta x, \forall i$ and $ \Delta y_i = \Delta y, \forall j$.  For the illustration, we assume that the following $3$ cell averages,  ${\bar v}_{i,j}$, ${\bar v}_{i+1,j}$ and ${\bar v}_{i,j-1}$ are determined to be used for the reconstruction by the Newton's divided difference method.  We only show the reconstruction at the cell boundary of  $(x,y)  = (x_{i+ {1\over 2}}, y_j)$. The reconstruction for other cell boundaries can be achieved similarly.

\subsection{Polynomial reconstruction}
Assume that the following interpolant $s(x,y)$ could generate the given $3$ cell averages ${\bar v}_{i,j}$, ${\bar v}_{i+1,j}$, ${\bar v}_{i,j-1}$
$$ s(x,y) = \lambda_1 + \lambda_2 x + \lambda_3 y. $$
Let $ {\vec{V}} =  [{\bar v}_{i,j}, {\bar v}_{i+1,j}, {\bar v}_{i,j-1}]^T$, $ {\vec{\lambda}} = [\lambda_1, \lambda_2,  \lambda_3]^T$ and the interpolation matrix $A$ be
$$
A = \begin{bmatrix} 
1 & x_i & y_j+dy\\ 
1 & x_i & y_j\\ 
1 & x_i+dx & y_j
\end{bmatrix}. 
$$
Then the expansion coefficients $\lambda_i$ are given by solving the linear system $ \vec{V} = A \cdot \vec{\lambda} $.
After plugging $\lambda_i$ in $s(x,y)$ at $(x,y) = (x_{i+ {1\over 2}}, y_j)$, we obtain
\begin{equation}
v^-_{i+ {1\over 2},j} = s(x_{i+ {1\over 2}},y_j) = {1\over 2} \cdot {\bar v}_{i,j}+{1\over 2} \cdot {\bar v}_{i+1,j} + 0 \cdot {\bar v}_{i,j-1}. 
\label{poly_recon2d}
\end{equation}
Expanding $v^-_{i+ {1\over 2},j}$ around $x = x_{i+ {1\over 2}}$ and $y = y_j$ in the Taylor series yields
\begin{equation}
v^-_{i+ {1\over 2},j}  = v(x_{i+ {1\over 2}},y_j) + {1\over 6}  v_{xx}(x_{i+ {1\over 2}},y_j) \Delta x^2 + {1\over 24} v_{yy}(x_{i+ {1\over 2}},y_j) \Delta x^2+ O(\Delta x^3)+O(\Delta y^3). 
\label{poly_recon_error2d}
\end{equation}
$ v(x_{i+ {1\over 2}},y_j) $ in \eqref{poly_recon_error2d} is the exact value of $v(x)$ at $x = x_{i+ {1\over 2}}$ and $y = y_j$. So we confirm that \eqref{poly_recon2d} is a 2nd order reconstruction.

\subsection{Perturbed polynomial reconstruction}
The 2D RBF reconstruction can be obtained as a straightforward extension of the 1D RBF reconstruction. However, the calculation can be complicated due to the double integrals of RBFs. For this reason, we consider this problem from a different perspective. That is, instead of directly constructing the 2D RBF interpolation with a specific RBF basis function, we modify the polynomial reconstruction by adding small perturbation terms containing the shape parameter as we did in Section \ref{perturbation}.  As in the 1D reconstruction, we then optimize the shape parameter to increase the order of convergence.

Mimicking the 1D RBF interpolation, we assume that the perturbed 2D polynomial interpolation, (\ref{poly_recon2d}), is given by the following form
\begin{equation}
v^-_{i+ {1\over 2},j}  = \left({1\over 2}+c_1 \epsilon^2 h^2\right) \cdot {\bar v}_{i,j}+\left({1\over 2}+c_2 \epsilon^2 h^2\right) \cdot {\bar v}_{i+1,j} + \left(0+c_3 \epsilon^2 h^2\right) \cdot {\bar v}_{i,j-1}, 
\label{rbf_recon_2d}
\end{equation}
where $h = h(\Delta x, \Delta y)$. Expanding $v^-_{i+ {1\over 2},j}$ around $x = x_{i+ {1\over 2}}$ and $y = y_j$ in the Taylor series yields
\begin{eqnarray}
v^-_{i+ {1\over 2},j} &=& v(x_{i+ {1\over 2}},y_j) + {1\over 6}  v_{xx}(x_{i+ {1\over 2}},y_j) \Delta x^2 + {1\over 24} v_{yy}(x_{i+ {1\over 2}},y_j) \Delta y^2\nonumber \\
&&+ (c_1+c_2+c_3) \epsilon^2 h^2 v(x_{i+ {1\over 2}},y_j)\nonumber \\
&&+  O(\Delta x^3)+O(\Delta y^3). 
\label{rbf_recon_error_2d}
\end{eqnarray}
Thus if we take the value of $\epsilon$ as below
\begin{equation}
\epsilon^2 = -\frac{{1\over 6}  v_{xx}(x_{i+ {1\over 2}},y_j) \Delta x^2 + {1\over 24} v_{yy}(x_{i+ {1\over 2}},y_j)\Delta y^2}{(c_1+c_2+c_3) v(x_{i+ {1\over 2}},y_j)h^2 },
\label{i+1/2_eps_2d}
\end{equation}
then we obtain a $3$rd order accurate approximation.
Again, the exact values of $v_{xx}, v_{yy}$ and $v(x_{i+{1\over 2}}, y_j)$ in (\ref{i+1/2_eps_2d}) are not available, so we replace them with their approximations based on the given cell averages to approximate $\epsilon^2$. Fortunately, it is possible to approximate all those quantities using the given $5$ cell averages. For example, $v_{xx}$ and $v_{yy}$ can be easily approximated using $\left\{ {\bar v}_{i-1,j}, {\bar v}_{i,j}, {\bar v}_{i+1,j}\right\}$ and $\left\{ {\bar v}_{i,j-1}, {\bar v}_{i,j}, {\bar v}_{i,j+1}\right\}$, respectively with a $2$nd order accuracy (see Figure \ref{figure2d_grid}). Then for this case, the $\epsilon^2$ is approximated by 
\begin{equation}
\epsilon^2 \approx  \frac{{\bar v}_{i,j+1}+4{\bar v}_{i-1,j}-10{\bar v}_{i,j}+4{\bar v}_{i+1,j}+{\bar v}_{i,j-1}}{4 (c_1+c_2+c_3)(-2 {\bar v}_{i-1,j}+{\bar v}_{i,j}-5{\bar v}_{i+1,j}) h^2 +\epsilon_M}.
\label{i+1/2_eps_approx_2d}
\end{equation}
To check whether this approximation of $\epsilon^2$ still achieves the $3$rd order accuracy, we expand $v^-_{i+ {1\over 2},j}$ around $x = x_{i+ {1\over 2}}$ and $y = y_j$ in the Taylor series. After a small calculation, we confirm that 
$$ v^-_{i+ {1\over 2},j} = s(x_{i+ {1\over 2}},y_j) = v(x_{i+ {1\over 2}},y_j) + O(\Delta x^3)+O(\Delta y^3). $$ 
Here note that for the $3$rd order approximation, we used up all the given cell averages. The polynomial interpolation, even in the case that all the cell averages are used, still yields a $2$nd order approximation because only $5$ cell averages are used. The perturbation terms in (\ref{rbf_recon_2d}) indeed yield the flexibility to use all the possible approximation from the given $5$ cell averages. 

\section{2D Numerical examples}
For the 2D numerical examples, we consider the 2D hyperbolic conservation laws
$$
    v_t + f(v)_x + g(v)_y = 0. 
$$
In each cell $I_{i,j}$, we have 
\begin{eqnarray}
    \int_{I_{i,j}} u_t(x,y,t)dx dy &=& 
     - 
     \int_{y_{j-{1\over 2}}}^{y_{j+{1\over 2}}} \left( f(u(x_{i+{1\over 2}}, y) - f(u(x_{i-{1\over 2}}, y) \right) dy \nonumber \\
     && - 
     \int_{x_{i-{1\over 2}}}^{x_{i+{1\over 2}}} \left( g(u(x,y_{j+{1\over 2}}) - g(u(x, y_{j-{1\over 2}}) \right) dx.  \label{2DIntegral}
\end{eqnarray}
Thus the 2D ENO/WENO finite volume method involves the approximation to integrals and the overall order provided by the method  depends not only on the value of $k$, but also on the number of the Gaussian quadrature points used for the integrals. If the method is only of $2$nd order accurate (polynomial interpolation with $k = 2$), then one quadrature point at the cell boundaries is enough for the integrals to maintain the same order. The RBF-ENO reconstruction is of the $3$rd order accurate, so we need to use at least two quadrature points for the desired order. If one quadrature point is used as the regular ENO method instead, however, the RBF-ENO still maintains the $2$nd order accuracy, but it yields higher accuracy than the regular ENO method because the RBF-ENO reconstruction is already $3$rd order accurate. Since we only want to modify the existing $2$nd order regular ENO code with the minimum changes, we still use one quadrature point and yet improve much the accuracy. 
\subsection{Numerical results}
First we check the reconstruction error by the RBF-ENO method applied to a smooth function $u(x) = \sin(2 \pi (x+y))$ with $k = 2$ to confirm the desired order of convergence. If five cell averages, ${\bar v}_{i-1,j}$, ${\bar v}_{i,j}$, ${\bar v}_{i+1,j}$, ${\bar v}_{i,j-1}$, and ${\bar v}_{i,j+1}$ are available, the regular ENO method chooses three of them by the Newton's divided difference method, which ends up with the $2$nd order convergence. The RBF-ENO method also chooses three cell averages from the Newton's divided difference, but the five cell averages are used to optimize the shape parameter. This will provides the $3$rd order accuracy for the smooth problem. For comparison, we also use the polynomial reconstruction using all the five cell averages. If we assume the uniform grid with $dx = dy$, the possible interpolation using the $5$ cells is given by 
\begin{eqnarray}
   s(x,y) = c_0 + c_1 x + c_2 y + c_3 x^2 + c_4 y^2, 
   \label{FVM_interpolation}
\end{eqnarray}
where $c_i, i = 0, \cdots,  4,$ are all constants determined by the given five cell averages. Here note that there is no cross term, $xy$, due to the axial symmetry. Thus we only expect at most the $2$nd order convergence although those $5$ cell averages are all used while the RBF-ENO can yield the $3$rd order convergence. 

\begin{figure}[h]
\begin{center}
\includegraphics[width=0.47\textwidth]{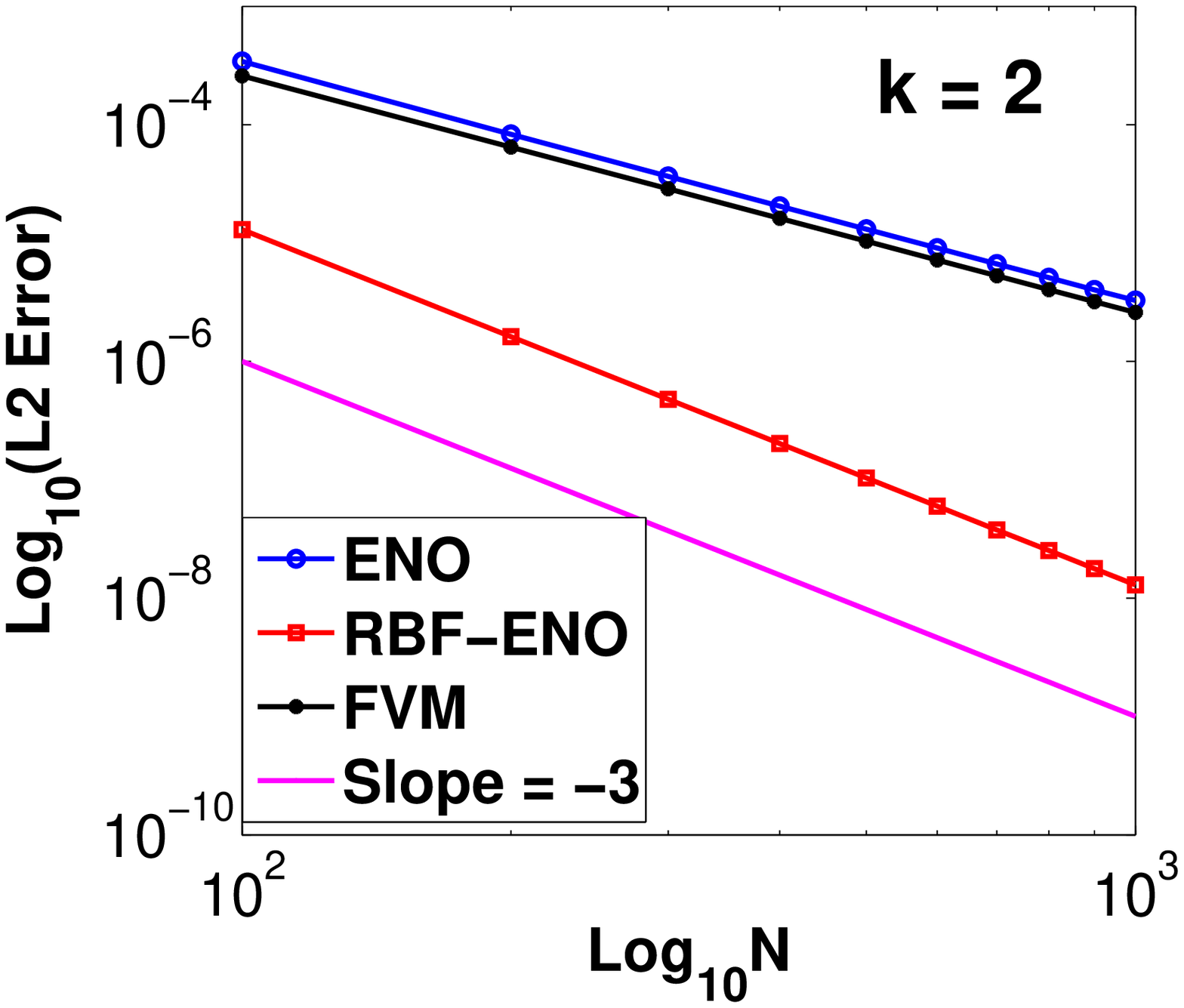}
\end{center}
\caption{(Color online).  $L_2$ errors versus $N$ on logarithmic scale with $k = 2$. ENO method (blue), RBF-ENO method (red), 5-point Finite Volume method (black).}
\label{test1_2d}
\end{figure}

Figure \ref{test1_2d} shows the $L_2$ errors  versus $N$ on logarithmic scale for the regular ENO (blue circle), RBF-ENO (red square), and the 5-cell finite volume method (black filled circle) reconstructions. For the 5-cell finite volume method, (\ref{FVM_interpolation}) is used. The solid line in magenta is the reference line of order $3$. As shown in the figure, the RBF-ENO reconstruction yields the desired $3$rd order convergence and is much more accurate than the regular ENO method or the 5-cell finite volume method. The regular ENO and the 5-cell finite volume methods yield only the $2$nd order convergence. We also note that the difference between the 5-cell finite volume reconstruction with (\ref{FVM_interpolation}) and the regular ENO reconstruction is not significant. 

%
%
%
%

\subsubsection{Example 1}
We solve the advection equation for $(x,y) \in [0,1]\times [0,1]$
\begin{eqnarray}
v_t + v_x + v_y = 0, \quad  t > 0,
\label{example_2d_1}
\end{eqnarray}
with the initial condition $u (x,0) = u_0(x) = \sin(2 \pi (x+y))$ and the periodic boundary condition. The CFL condition is given by $\Delta t \le C \max(\Delta x, \Delta y)$ with $C = 0.1$. 

Table \ref{table0_2d} shows the $L_1$, $L_2$, and $L_\infty$ errors at the final time, $T = 0.5$ for the regular ENO method, 5-cell finite volume, and the RBF-ENO method. For the integral in the RHS of (\ref{2DIntegral}), we used one quadrature point at the mid point of each cell boundary, which makes all those methods yield the $2$nd order convergence at most. From the table, we observe that the RBF-ENO is more accurate than the regular ENO or the 5-cell finite volume method. This is again because the reconstruction at the quadrature points by the RBF-ENO is $3$rd order accurate. 

\begin{table}[h]
\renewcommand{\arraystretch}{3.5}
\caption{The $L_1, L_2$ and $L_\infty$ errors for the linear advection equation, \eqref{example_2d_1} with $k = 2$ at $T = 0.5$. }
\begin{center} \footnotesize
\renewcommand{\arraystretch}{1.1}
\begin{tabular}{|c|c|c|c|c|c|c|c|} 
\hline  
Method & N & $L_1$ error & $L_1$ order & $L_2$ error & $L_2$ order & $L_\infty$ error & $L_\infty$ order\\ 
\hline 
       & 10 & 9.90E-2 &   --   & 1.10E-1 &    --  & 1.74E-1 & --     \\  
       & 20 & 3.86E-2 & 1.3586 & 4.26E-2 & 1.3697 & 7.99E-2 & 1.1230 \\  
ENO    & 40 & 1.08E-2 & 1.8351 & 1.40E-2 & 1.6026 & 3.43E-2 & 1.2164 \\
       & 80 & 3.00E-3 & 1.8492 & 4.47E-3 & 1.6526 & 1.43E-2 & 1.2594 \\ 
       & 160& 8.01E-4 & 1.9078 & 1.42E-3 & 1.6539 & 5.91E-3 & 1.2809 \\ 
       & 320& 2.12E-4 & 1.9177 & 4.52E-4 & 1.6517 & 2.40E-3 & 1.2963 \\  
\hline
       & 10 & 5.86E-2 &   --   & 5.94E-2 &    --  & 7.53E-2 & --     \\  
       & 20 & 1.43E-2 & 2.0338 & 1.51E-2 & 1.9750 & 2.03E-2 & 1.8917 \\  
5 point& 40 & 3.48E-3 & 2.0369 & 3.77E-3 & 2.0016 & 5.19E-3 & 1.9665 \\
finite volume    & 80 & 8.58E-4 & 2.0232 & 9.42E-4 & 2.0034 & 1.31E-3 & 1.9817 \\ 
       & 160& 2.13E-4 & 2.0100 & 2.35E-4 & 2.0021 & 3.30E-4 & 1.9933 \\ 
       & 320& 5.30E-5 & 2.0055 & 5.87E-5 & 2.0012 & 8.28E-5 & 1.9967 \\  
\hline
           & 10 & 1.31E-2 &    --  & 1.41E-2 &    --  & 2.20E-2 &    --  \\  
           & 20 & 3.52E-3 & 1.8986 & 3.77E-3 & 1.9013 & 5.70E-3 & 1.9493 \\  
RBF-ENO    & 40 & 8.57E-4 & 2.0411 & 9.40E-4 & 2.0060 & 1.30E-3 & 2.1301 \\
           & 80 & 2.10E-4 & 2.0266 & 2.31E-5 & 2.0207 & 3.21E-4 & 2.0203 \\ 
           & 160& 5.20E-5 & 2.0159 & 5.74E-6 & 2.0108 & 8.08E-5 & 1.9900 \\ 
           & 320& 1.29E-5 & 2.0086 & 1.43E-7 & 2.0050 & 2.02E-5 & 1.9995 \\
\hline
\end{tabular}
\end{center} 
\label{table0_2d}
\end{table}

\subsubsection{Example 2}
\begin{figure}[h]
\begin{center}
\includegraphics[width=0.7\textwidth]{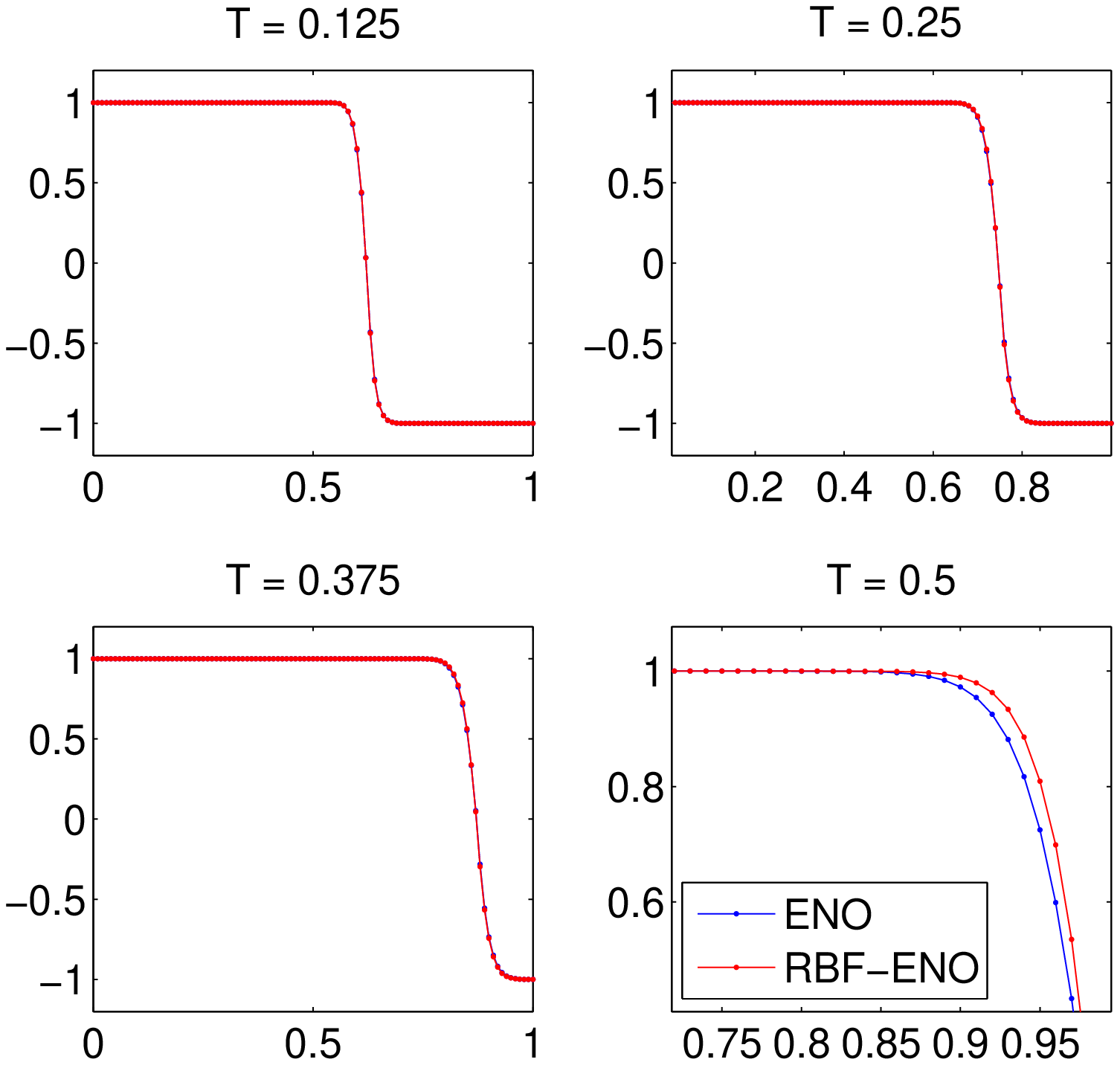}
\end{center}
\caption{(Color online). Solutions to \eqref{example_2d_1} at $x=0.5$ with the discontinuous initial condition, \eqref{example_2d_2} for the ENO (blue), RBF-ENO (red) methods with $k = 2$ and $N = 100$.}
\label{figure_2d_1}
\end{figure}

We consider the same linear advection equation \eqref{example_2d_1} with the following discontinuous initial condition and the boundary condition 
\begin{eqnarray}
v(x,y,0) &=& \left\{
             \begin{array}{rcl}
             1 &\text{if} &y \in [0,0.5] \\
             -1 &\text{if} &y \in (0.5,1]  
             \end{array}  
        \right. ,\\
        v(x,0,t) &=& v(x,1,t) = 1, \quad  t>0. 
\label{example_2d_2}        
\end{eqnarray}
With this example, we check how the RBF-ENO solution behaves near the discontinuity. To handle the discontinuous solution, we apply the same monotone polynomial method as in 1D. In Figure \ref{figure_2d_1}, the solutions by the regular ENO (blue) and RBF-ENO (red) methods are given with time at $x = 0.5$. The figure clearly shows that the RBF-ENO yields non-oscillatory solutions for all time. The figure also shows that the RBF-ENO method yields much sharper solution profile than the regular ENO method, which is clearly shown in the zoomed profile in the bottom right figure. 


\subsubsection{Example 3}
\begin{figure}[h]
\begin{center}
\includegraphics[width=0.5\textwidth]{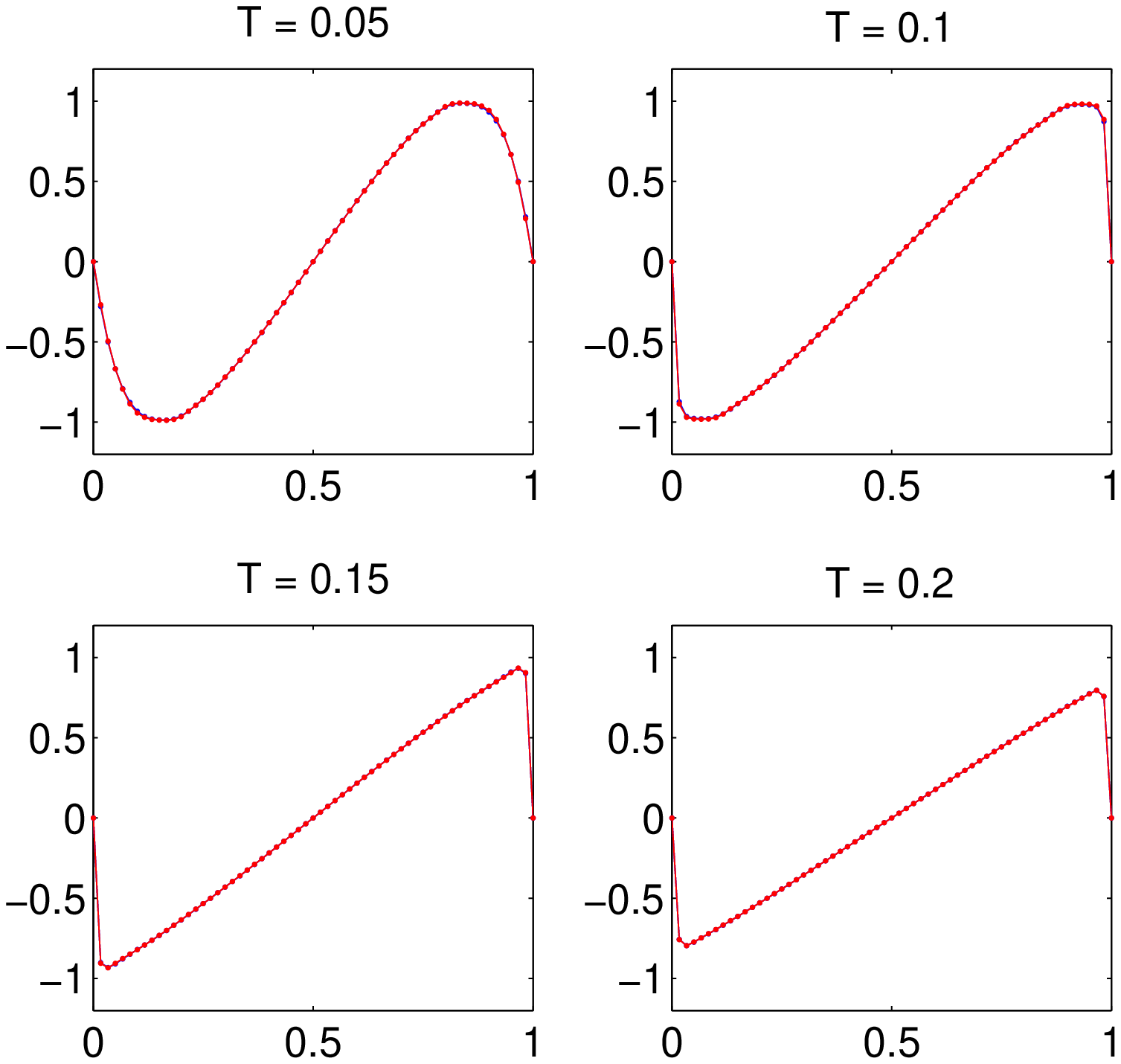}
\includegraphics[width=0.4\textwidth]{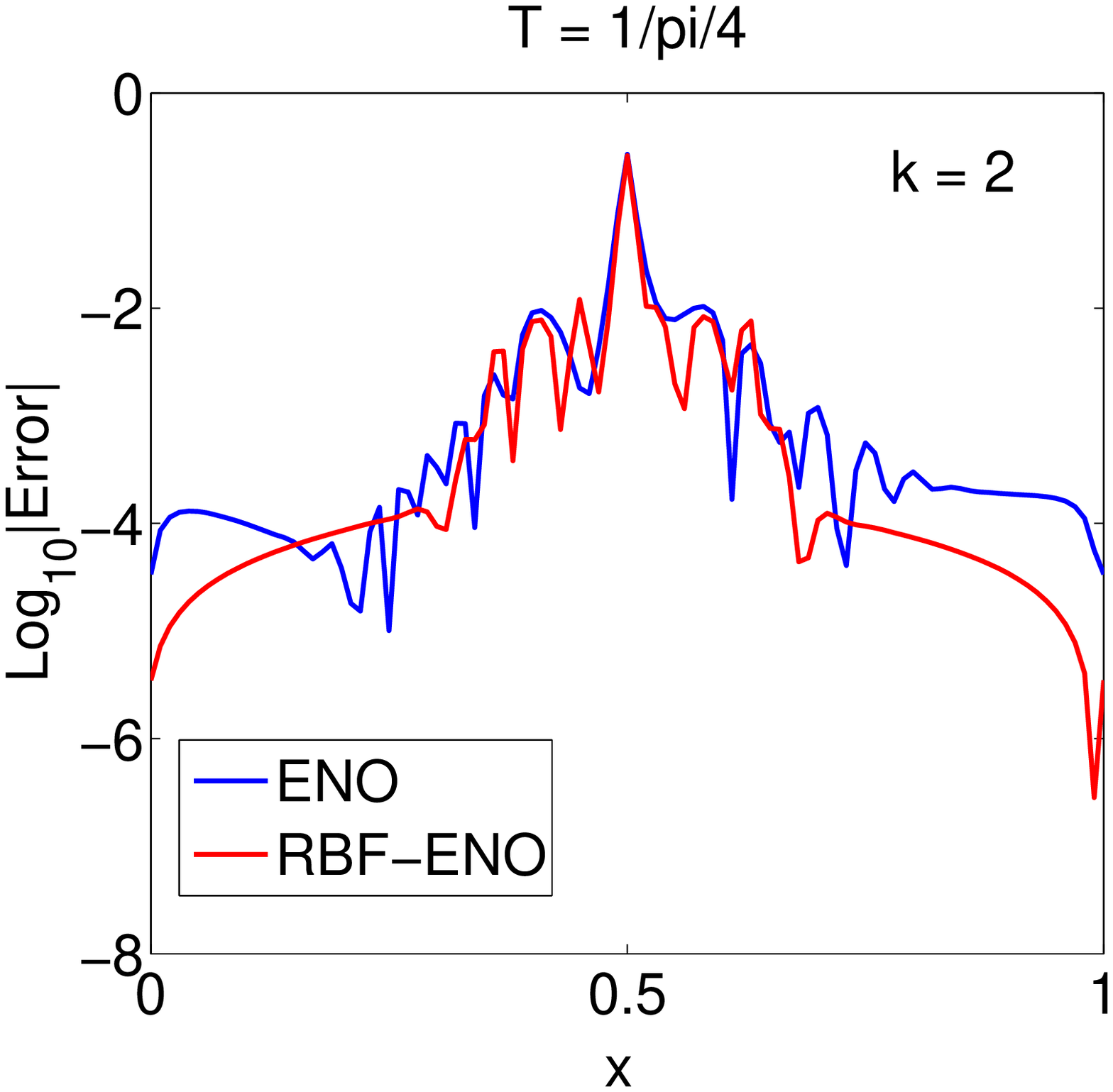}
\end{center}
\caption{(Color online). Left: The RBF-ENO solutions to \eqref{example3_2d} at various time. Right: The pointwise errors at $t = 1/4\pi$ by the ENO (blue) and RBF-ENO (red) on logarithmic scale with $k = 2$ and $N = 100$. }
\label{2d_figure3}
\end{figure}

We consider the 2D Burgers' equation for $(x,y) \in [0,1]\times [0,1]$ 
\begin{eqnarray}
   v_t + \left(\frac{1}{2} v^2 \right)_x +\left(\frac{1}{2} v^2 \right)_y = 0,  \quad t >0, 
\label{example3_2d}
\end{eqnarray}
with the initial condition 
\begin{eqnarray}
        v(x,0) = \sin(2 \pi (x+y) ). 
\label{example3_example_2d}
\end{eqnarray}
Figure \ref{2d_figure3} shows the RBF-ENO solutions at various times at $x = 0.5$ (left figure) and the pointwise errors by the ENO (blue) and RBF-ENO (red) methods with $k = 2$ and $N = 100$ at $t = 1/4\pi$. The left figure clearly shows that the RBF-ENO solution is not oscillatory yet yielding a sharp shock profile near the boundaries. The right figure shows that the RBF-ENO method yields more accurate results than the regular ENO method in the smooth region. 
%
%

\section{Conclusion}
\label{sec:conclusion}
In this paper, we developed a non-polynomial ENO method for solving hyperbolic equations.  As an example of non-polynomial bases, we used RBFs. The formulation based on the non-polynomial basis yields the flexibility of improving the original ENO accuracy. The key idea of the developed method lies in the adaptation of the shape parameters in the expansion with a non-polynomial basis that can make the leading error term vanish or at least become small in the local interpolation. The new non-polynomial ENO method improves local accuracy and convergence if the underlying solution is smooth. For the non-smooth solution such as a shock, we adopted the monotone interpolation method so that the non-polynomial ENO reconstruction is reduced into the regular ENO reconstruction resulting in the suppression of the Gibbs oscillations. 
The numerical results show that the non-polynomial ENO method is superior to the regular ENO method and even better than the WENO-JS method for $k = 2$. The non-polynomial ENO method yields  $4$th order accuracy while the regular ENO solution is only $3$rd order accurate for $k = 3$ in the smooth region. The numerical results show that the developed method yields highly accurate results for the scalar problems for both $k = 2$ and $k = 3$. For the system problems, the non-polynomial ENO solutions are similar to or better than the regular ENO solutions. The non-polynomial ENO scheme is slightly more costly than the regular ENO scheme because it has a procedure of computing the optimal shape parameter values. But it is less expensive than the WENO scheme for the given value of $k$. For some cases, the non-polynomial ENO method achieves the same level of accuracy as the WENO method or even better accuracy than the WENO method while its computational cost is less demanding than the WENO method. The WENO method based on the non-polynomial ENO reconstruction is also better then the regular WENO method. The 2D non-polynomial finite volume interpolation is more beneficial than the 2D polynomial interpolation. We showed that the non-polynomial interpolation can raise the order which can not be obtained with the polynomial interpolation even though all the cell averages are used. We provided the table of the reconstruction coefficients for $k = 2$ and $k = 3$ for the non-polynomial ENO method. The non-polynomial ENO formulation for higher values of $k$ will be considered in our future work. The current work considered the uniform mesh only. In our future work, we will investigate the non-polynomial ENO method with the nonuniform and unstructured mesh. As mentioned in Introduction, the meshless feature of RBFs was combined with the WENO method in \cite{RBF-ENO} where the shape parameter was globally fixed for the reconstruction. It will be interesting to investigate how the optimization of the shape parameter can be realized with the meshless properties of RBFs on the unstructured mesh. 



\vskip .1in
\noindent
{\bf Acknowledgments:} The authors thank W.-S. Don for his useful comments on the construction of the RBF-ENO/WENO method. The second author thanks Grady Wright for his useful comments on the RBF interpolation.  


\end{document}